\documentclass[final,leqno,onefignum,onetabnum]{siamltex}
\usepackage{fullpage}
\usepackage{comment}
\usepackage{upgreek}
\usepackage[mathscr]{eucal}
\usepackage{color}
\usepackage[usenames,dvipsnames,svgnames,table]{xcolor}
\usepackage{amsmath,amssymb,amsopn,mathtools}
\usepackage{graphicx}
\usepackage{float}
\usepackage{adjustbox}
\usepackage{diagbox}
\usepackage{multirow}
\usepackage{makecell}
\usepackage{hyperref}       
\hypersetup{
colorlinks=true,   
citecolor=	PineGreen,     
	filecolor=PineGreen,     
	linkcolor=PineGreen,    
	urlcolor= PineGreen}      

\usepackage{relsize}
\numberwithin{equation}{section}
\usepackage{enumitem}

   \allowdisplaybreaks
   \DeclareMathOperator*{\argmin}{argmin}

\newcommand{\tu}{\textup}

\usepackage{slashbox}
\usepackage{lipsum}
\newcommand{\bfa}[1]{\boldsymbol{#1}} 			%

\usepackage{color}
\usepackage{adjustbox}
\usepackage{diagbox}
\definecolor{black}{rgb}{0,0,0}

\definecolor{red}{rgb}{1,0,0}

\definecolor{blue}{rgb}{0,0,1}

\usepackage{subcaption}
\renewcommand{\theequation}{\arabic{section}.\arabic{equation}}
\numberwithin{equation}{section}

\renewcommand{\div}{\mathop{\rm div}\nolimits}

\newcommand{\todo}[1]{{\color{red}{#1}}}

\newcommand{\ds}[1]{{\color{blue}{#1}}}

\newcommand{\dx}{ \mathrm{d}x}

\newcommand{\dt}{ \mathrm{d}t}

\newcommand{\beq}{\begin{equation}}
\newcommand{\eeq}{\end{equation}}
\newcommand{\beqq}{\begin{equation*}}
\newcommand{\eeqq}{\end{equation*}}
\newcommand{\beqas}{\begin{eqnarray*}}
\newcommand{\eeqas}{\end{eqnarray*}}
\newcommand{\bsp}{\begin{split}}
\newcommand{\eesp}{\end{split}}

\usepackage[english]{babel}

\usepackage{tikz}

\usetikzlibrary{shapes.geometric, arrows}
\usetikzlibrary{positioning}
\usetikzlibrary{calc}

\tikzstyle{Input} = [rectangle, minimum height=2.2cm, text width=4.5cm, text centered, draw=black, fill=red!30]
\tikzstyle{DNN} = [rectangle, rounded corners, minimum width=3cm, minimum height=2.2cm, text width=4.5cm, text centered, draw=black, fill=blue!30]
\tikzstyle{MS} = [rectangle, minimum width=3cm, minimum height=2.2cm, text width=4.5cm, text centered, draw=black, fill=orange!30]

\tikzstyle{arrow} = [thick,->,>=stealth]

\makeatletter
\def\thanks#1{\protected@xdef\@thanks{\@thanks
		\protect\footnotetext{#1}}}
\makeatother

\begin{document}
\title{ Prediction of discretization of online GMsFEM using deep
learning for Richards equation}

\author{Denis Spiridonov$^{\MakeLowercase{a}}$\thanks{\textit{Denis Spiridonov}; $^a$Laboratory of Computational Technologies for Modeling Multiphysical and Multiscale Permafrost Processes, North-Eastern Federal University, 677000 Yakutsk, Republic of Sakha (Yakutia), Russia; \texttt{d.stalnov@mail.ru}
}
\and
Sergei Stepanov$^{\MakeLowercase{a}}$\thanks{\textit{Sergei Stepanov}; $^a$Laboratory of Computational Technologies for Modeling Multiphysical and Multiscale Permafrost Processes, North-Eastern Federal University, 677000 Yakutsk, Republic of Sakha (Yakutia), Russia; \texttt{cepe2a@inbox.ru}
}
\and
Tina Mai$^{\MakeLowercase{b,c,*}}$\thanks{$^{*}$Corresponding author: \textit{Tina Mai}; $^b$Institute of Research and Development, Duy Tan University, Da Nang, 550000, Vietnam; $^c$Faculty of Natural Sciences, Duy Tan University, Da Nang, 550000, Vietnam; \texttt{maitina@duytan.edu.vn}
}
}
\maketitle

\maketitle

\begin{abstract}
We develop a new coarse-scale approximation strategy for the nonlinear single-continuum Richards equation as an unsaturated flow over heterogeneous non-periodic media, using the online generalized multiscale finite element method (online GMsFEM) together with deep learning.
A novelty of this approach is that local online multiscale basis functions are computed rapidly and frequently by utilizing deep neural networks (DNNs).  More precisely, we employ the training set of stochastic permeability realizations and the computed relating online multiscale basis functions to train neural networks. The nonlinear map between such permeability fields and online multiscale basis functions is developed by our proposed deep learning algorithm.  That is, in a new way, the predicted online multiscale basis functions incorporate the nonlinearity treatment of the Richards equation and refect any time-dependent changes in the problem's properties.  Multiple numerical experiments in two-dimensional model problems show the good performance of this technique, in terms of predictions of the online multiscale basis functions and thus finding solutions.
\end{abstract}



\section{Introduction}\label{intro}
In general soils, unsaturated flow occurs when water seeps through many soil pores, evaporation dries up the topsoil, or water moves downward \cite{unflow, soilmoisture}. This typically results in the idea of soil moisture \cite{soilmoisture}, which is defined by water potential (the energy state of water) and water content (the volume of water held in the spaces between unsaturated soil particles for any given soil sample). Soil moisture is essential for several hydrological, ecological, and biogeochemical procedures, even though it is not a major part of the water cycle. That is to say, weather forecasting, agronomy, sanitary engineering, environmental management,  groundwater storage, energy balances, earth dynamic systems, etc.\ depend heavily on soil moisture.  Richards equation \cite{r0,richarde1,richarde2,haverkamp,richardsreview} describes mathematically the related processes, showing how water flows through unsaturated soil surface that is filled with both air and water \cite{ry1, soilmoisture}.

Given non-periodic and heterogeneous media, the Richards equation is examined in this study.  As with subsurface operations, oil reservoirs, aquifers, infiltration, and most real-world simulation problems, obstacles exacerbate when there are major property uncertainties in some local locations.  We emphasize that such uncertainties usually vary spatially.

These kinds of challenges need extremely fine grids, where the popular frameworks are the Finite Volume Method (FVM) and the Finite Element Method (FEM) \cite{fem}.  In our situation, such methods produce big systems with a lot of unknowns, which are computationally difficult and costly to solve. Consequently, some form of model reduction is necessary. That is, typical approaches include upscaling based on homogenization \cite{homo1,homod}, whose expansion can be used to derive multi-continuum methods in recents papers \cite{yewl22homo,eyjw2023mchomo,dsdw2023mchomo}.  Also, as a part of our previous work \cite{sdt22}, numerical homogenization was helpful to handle the Richards equation.

In our current study, scale separation is not required.  Thus, avoiding classical homogenization process, we consider the generalized multiscale finite element method (GMsFEM) to propose a good numerical strategy for solving the stochastic heterogeneous nonlinear Richards equation.  As an overview, creating local coarse-grid multiscale basis functions is the main objective of the GMsFEM. To achieve this, local snapshot spaces (comprising a few local solutions) are built, and inside these snapshot spaces, local spectral decomposition is performed. Hence, the generated eigenfunctions are able to transmit the local features to the global ones via those coarse-grid multiscale basis functions.
Next, by solving a coarse-grid matrix problem in the multiscale space (created by such multiscale basis functions), the approximations of solution are reached.

For the Richards equation, the GMsFEM was successfully used in the recent articles of one of the authors: the work with offline multiscale basis functions \cite{Spiridonov2019}, and the work utilizing online multiscale basis functions \cite{denisonr}.
More specifically, in the offline GMsFEM (as in \cite{Spiridonov2019}), multiscale basis functions (offline function space) are generated 
and then employed to solve the Richards equation on a coarse grid. If the multiscale basis functions depend on online data (such as the parameters or right-hand side), then they form online function spaces \cite{yeb2}.  The latter approach, known as online GMsFEM (as in \cite{denisonr}), follows the residual-driven online GMsFEM \cite{chungres}, where the establishment of online multiscale basis functions is based on local residuals in relation to adaptivity \cite{chungres,gne}, which refers to the insertion of online multiscale basis functions into specific locations. In our paper, instead of adaptivity, it is sufficient to uniformly add one online multiscale basis function to each coarse neighborhood.  The goal of building such online multiscale basis functions is to locally minimize errors. 
 Theoretical frameworks demonstrate that these functions provide quick convergence \cite{chungres}.

Moreover, for nonlinear problems, the online GMsFEM is an appropriate choice, as it introduces additional online multiscale basis functions to include any time-dependent changes in the problem's properties \cite{chungres, online29nl, denisonr}. The online GMsFEM showed good performance in one of the authors' previous work \cite{denisonr}, where one can observe superiority of the online GMsFEM over the traditional offline GMsFEM, in solving the Richards equation describing unsaturated infiltration.  Within our paper, for nonlinear single-continuum Richards equation, we combine the online GMsFEM with the idea of Picard iterative process \cite{gne}.  Specifically, the offline GMsFEM will be utilized throughout each Picard iteration of linearization, then with these data, the online GMsFEM will need to be employed only at the final Picard step (using the online multiscale basis functions obtained from the previous Picard step).
In this time-saving manner, the online GMsFEM can address the issues posed by the media's multiple scales, heterogeneities, and high contrast. 

However, more challenges occur when there are uncertainties regarding the media characteristics in specific local areas, which is typical for aquifers or oil reservoirs.  Sampling realizations of media features is a simple way to quantify such uncertainties. With respect to each permeability field, a corresponding local online multiscale basis function is created by the online GMsFEM (see Section \ref{onstage}).  When there are many different realizations of permeability fields, to obtain these basis functions, the computational load can become tremendous. Thus, creating a useful connection between the permeability fields and the local online multiscale basis functions can greatly reduce the computational complexity, by avoiding the need for expensive, time-consuming calculations.  A functional relation of this kind is typically nonlinear due to different medium configurations, leading to high-order approximations during the modeling process.  Thus, it makes sense to use machine learning methods to build and solve these intricate models \cite{dnngms19, ex-imML,sdt22}.  

One kind of machine learning frameworks is formed by deep neural networks (DNNs) \cite{dnneg14}. In this study, local online multiscale basis functions (via the online GMsFEM) are examined using DNNs. The architecture of a deep neural network is normally composed of multiple layers, with many neurons (as nonlinear processing units for property extraction) in each layer \cite{dnnnature}. Following the input layer (usually without neurons), a deep neural network (DNN)'s first hidden layer learns fundamental traits that are then sent to the second hidden layer so that it can train itself to recognize more abstract and complicated traits. This structure is continued on next levels as the layers are added until the output is sufficiently precise to be deemed admissible.  In order to improve the neural networks' expressiveness, certain nonlinear activation functions (like ReLU, tanh, and sigmoid) are needed between layers, from an input neuron to an output neuron.  After that, the output is transferred as an input to the network's subsequent layer.  DNNs have demonstrated efficacy in solving pattern recognition tasks, such as speech recognition, picture identification, and natural language operation; and DNNs have also been employed to understand sophisticated data sets \cite{dnneg8,dnneg9,dnneg10}.  To further explore the depictive capabilities of deep neural networks, considerable research has been accomplished \cite{dnneg11,dnneg12,dnneg13,dnneg14,dnneg15}. This work is also inspired by numerous other uses of convolutional neural networks (CNNs) \cite{agss,cnn-homo-elasticity19,cnn-homo-poro20} and DNNs \cite{nlnlmc31,wang2020deep,wang2020reduced, dnneg20}, especially our recent article \cite{sdt22}.

The main motivation for our utilization of deep neural networks is to speed up the procedure of computing local online multiscale basis functions. In the online GMsFEM, calculating offline multiscale basis functions needs to be carried out only once at the offline stage. Whereas, the creating of online multiscale basis functions must be performed many times at the online stage. This process can take a lot of execution time, so we replace this ``slow'' part of the online GMsFEM with a deep neural network.

That is to say, as an originality of our paper, we apply neural networks to the prediction of local online multiscale basis functions for a fast and frequent protocol.
Here, a deep neural network and graphics processing unit (GPU) are employed to construct a machine learning approach through several two-dimensional experiments, which are conducted during the training process. As a result, the predicted local online multiscale basis functions provide a new insight into the treatment of the nonlinearity of the Richards equation.  A matching set of targets (correct outputs) is used to train the neural networks with different stochastic permeability fields (inputs). Neural networks that have been trained on training data are capable of performing new calculations quickly and precisely.  Through various numerical tests, our results demonstrate that the established deep networks can be effectively applied to the testing data and yield good results. 

The structure of this paper is as follows.   In Section \ref{sec:model}, we begin with some preliminaries and then introduce the model problem for the nonlinear single-continuum Richards equation, that describes the unsaturated flow in the context of non-periodic heterogeneous media.  Section \ref{fineapp} discusses fine-scale discretization and Picard iteration for linearization.  Finding the coarse-grid approximation of the Richards equation's solution is the focus of Section \ref{coarseapp}, where both the offline and online GMsFEM are presented. Next, we introduce the idea of utilizing deep learning as a tool for predicting local online multiscale basis functions (instead of obtaining them by the online GMsFEM) through Section \ref{dlgms}, where the sampling is covered in detail, and the networks are precisely established. In Section \ref{numer}, to illustrate the effectiveness of our proposed networks via the prediction of online multiscale basis functions, we provide numerical examples (particularly for solutions), across multiple situations with different setups. The conclusion and summary are given in Section \ref{sec:conclusions}.  



\section{Model problem}\label{sec:model}
Let $\Omega$ be a simply connected, convex, open, bounded, Lipschitz computational domain in $\mathbb{R}^d\,.$  To clarify our reasoning throughout this work, we only examine the situation $d = 2\,;$ however, the idea is also easily extended to $d = 3\,.$  We use both Latin indices (as $i\,,j$) and Greek indices (like $\alpha\,,\iota$).  The spatial gradient and time derivative are respectively denoted by the symbols $\nabla$ and $\dfrac{\partial}{\partial t}\,.$  Additional notation can be found in \cite{sdt22,mcl,rtt21}. Over $\Omega\,,$ italic capitals (e.g.\ $f$) are employed to indicate functions, while bold letters (e.g.\ $\bfa{v}$ and $\bfa{T}$) denote vector fields and matrix fields. Also, italic capitals (e.g.\ $L^2(\Omega)$), boldface Roman capitals (e.g.\ $\bfa{V}$), and special Roman capitals (e.g.\ $\mathbb{S}$) represent the spaces of functions, vector fields, and matrix fields over $\Omega\,,$ respectively.

We denote by $(\cdot,\cdot)$ the $L^2$ inner product. The Sobolev space is symbolized by 
$V: = H_0^1(\Omega) = W_0^{1,2}(\Omega)= \overline{C_c^{\infty}(\Omega)}^{H^1(\Omega)}\,,$ where $C_c^{\infty}(\Omega)$ stands for the space of infinitely differentiable functions with compact support in $\Omega$ (see \cite{evans}, p.~256 and p.~259). The natural norm $\| \cdot \|_{V}$ of this space $V$ is as follows:
\[\|v\|_{V} = \left(\|v\|^2_{L^2(\Omega)} + 
\|\nabla v \|^2_{\bfa{L}^2(\Omega)}\right)^{\frac{1}{2}}\,.\] 
In this case, $\| \nabla v \|_{\bfa{L}^2(\Omega)}:= \| | \nabla v | \|_{L^2(\Omega)}\,,$ and the Euclidean norm of the two-component vector-valued function $\nabla v$ is indicated by $|\nabla v|\,.$  Additionally, let $\bfa{V} = V^2 = [H_0^1(\Omega)]^2\,.$  For every $\bfa{v} = (v_1,v_2) \in \bfa{V}\,,$ it is maintained that $\| \nabla \bfa{v}\|_{\mathbb{L}^2(\Omega)}:= \| | \nabla \bfa{v}| \|_{L^2(\Omega)}\,,$ with $|\bfa{A}|=\sqrt{\bfa{A} : \bfa{A}}$ serving as the Frobenius norm of the 
square matrix $\bfa{A}\,,$ having the Frobenius inner product defined as $\langle \bfa{A},\bfa{B}\rangle_{\tu{F}} = \bfa{A} : \bfa{B} = \tu{Tr}(\bfa{A}^{\tu{T}}\bfa{B}) = \sum_{i,j} a_{ij} b_{ij}\,.$ 

The norms of the Bochner space $L^r(0,T;X)$ \cite{evans} for each $1 \leq r < \infty$ are represented as
\begin{align}\label{bochner}
\begin{split}
\|\phi\|_{L^r(0,T;X)} &:= \left(\int_0^T \|\phi(t) \|_{X}^r \, \dt\right)^{1/r} < + \infty\,, \\
\|\phi\|_{L^{\infty}(0,T;X)} &: = \sup_{0 \leq t \leq T} \|\phi (t)\|_{X}  < + \infty\,.
\end{split}
\end{align}
At this point, $(X, \| \cdot \|_{X})$ is a Banach space; for instance, $X=H_0^1(\Omega) \,.$ We also define 
\[H^1(0,T;X):= \left \{ \phi \in L^2(0,T;X) \, : \, \partial_t \phi \in L^2(0,T;X) \right \}\,.\]
In order to simplify notation \cite{lporo}, the space for the pressure head $p(t,\cdot)$ is indicated by $V=H_0^1(\Omega)$ (as Bochner space) rather than $L^2(0,T;H_0^1(\Omega))\,.$

The time-dependent Richards equation \cite{r0} is typically expressed in the literature as
\beq
\label{eq:originaltheta}
\frac{\partial \Theta(p(t,\bfa{x}))}{\partial t} - \div (\varkappa(\bfa{x},p(t,\bfa{x}))\nabla p(t,\bfa{x}))
	= f(t,\bfa{x}) \  \textrm{in} \ (0,T] \times \Omega\,,
	\eeq
in which the terminal time is indicated by $T>0\,.$  We equip this equation with the Dirichlet boundary condition $p(t,\bfa{x})=0$ on $(0,T] \times \partial \Omega$ and the initial condition $p(0,\bfa{x})= p_{0}$ in $\Omega\,.$  The basic notation is available in \cite{rtt21}. Herein, $\Theta(p(t,\bfa{x}))$ represents the volumetric soil water content,
$p:=p(t,\bfa{x})$ stands for the pressure head, $\varkappa(p):=\varkappa(\bfa{x},p)$ shows the unsaturated hydraulic conductivity, and $f$  denotes the source or sink function.

Similar to \cite{rtt21}, we note that the volumetric water content $\Theta(p)$ of \eqref{eq:original0} usually takes the form of a nonlinear function of the pressure head $p$ as the following \cite{richarde1}:
\[\frac{\partial \Theta(p)}{\partial t} = C(p) \frac{\partial p}{\partial t}\,.\] 
For the online generalized multiscale finite element method (online GMsFEM) used in this research, the nonlinear hydraulic conductivity $\varkappa(\bfa{x}, p)$ is significant.  Contrarily, the component $C(p)$ can be removed as it is unnecessary, leaving us the identity function $\Theta(p) = p$ when we consider the Richards equation:
\beq
\label{eq:original0}
\frac{\partial p(t,\bfa{x})}{\partial t} - \div (\varkappa(\bfa{x},p(t,\bfa{x}))\nabla p(t,\bfa{x}))
= f(t,\bfa{x}) \  \textrm{in} \ (0,T] \times \Omega\,.
\eeq
Afterward, we assume that the multiscale high-contrast coefficient $\kappa(\bfa{x}) = \kappa(\bfa{x},\omega)$ in Eq.\ \eqref{eq:original0} is stochastic \cite{cnn-mcmc20}.  However, the stochastic symbol $\omega$ is eliminated for the sake of simplicity.

Since the hydraulic conductivity and its spatial gradient are assumed to be uniformly bounded, the following inequalities are satisfied by some positive constants $\underline{\varkappa}$ and $\overline{\varkappa}\,:$ 
\begin{align}
	\label{Coercivity}
	\begin{split}
		\underline{\varkappa} \leq \varkappa(\bfa{x}, p), \ | \nabla \varkappa(\bfa{x},p)| \leq \overline{\varkappa}\,.
	\end{split}
\end{align}
It is further assumed, without loss of generality, that the initial condition is \begin{equation}\label{ini} p_0 = p(0,\bfa{x}) \in V\,. \end{equation}

Using a given $u \in V\,,$ with any $p,v \in V\,,$ we define the bilinear form  
\begin{align}\label{ai}
a(p,v;u)=\int_{\Omega} \varkappa(u)\nabla p \cdot \nabla v \, \dx\,.
\end{align}

Subsequently, \eqref{eq:original0} has the variational form in the following manner: 
seek $p \in V$ such that 
\begin{align}\label{r1e}
	\left(\frac{\partial{p}}{\partial t} , v \right) + a(p,v;p) = (f,v)\,,
\end{align}
having all $v \in V\,,$ with a.e.\ $t \in (0,T]\,,$ and $f(t,\cdot) \in L^2(\Omega)\,.$ The initial condition was provided in \eqref{ini}.

\section{Fine-scale approximation and Picard iteration for linearization}
\label{fineapp}
To address the nonlinearity of our problem, we apply an efficient Picard iterative scheme as detailed in \cite{rpicardc,Spiridonov2019, gne, cemnlporo,tfcmm,rtt21,ttr22}. Such an iteration technique is given in this section for the time-dependent Richards equation.

The conventional backward Euler finite difference approach will be utilized to accomplish the first goal of \eqref{r1e}'s time discretization (see \cite{rpicardc,Spiridonov2019}, for instance): provided $p_s \in V\,,$ we seek $p_{s+1} \in V$ such that for every $v\in V\,,$
\begin{align}\label{r1ed}
	\left(\frac{p_{s+1} - p_{s}}{\tau} , v \right) + a(p_{s+1},v;p_{s+1})= (f_{s+1},v)\,. 
\end{align}
Here, the temporal range $[0,T]$ is equally divided into $S$ intervals, the size of the temporal step is $\tau = T/S > 0\,;$  the subscript $s$ (with $s=0,1,\ldots,S$) designates the evaluation of a function at the time $t_s= s\tau\,.$ 

The Picard linearization iteration (see \cite{rpicardc,Spiridonov2019,cemnlporo} for instance) will then handle the spatial nonlinearity as follows. Given $p_s \in V\,,$ guessing $p^0_{s+1} \in V$ at the $(s+1)$th time step. For $n=0,1,2, \ldots\,,$ find  
$p^{n+1}_{s+1} \in V$ such that for any $v \in V\,,$ 
\begin{align}\label{r1el}
\left(\frac{p^{n+1}_{s+1} - p_{s}}{\tau} , v \right) + a(p^{n+1}_{s+1},v;p^n_{s+1}) 
	= (f_{s+1},v)\,. 
\end{align}
This Eq.\ \eqref{r1el} is equivalent to
\begin{equation}\label{r1elr}
\left(\frac{p^{n+1}_{s+1}}{\tau}\,, v \right)  + a(p^{n+1}_{s+1},v;p^{n}_{s+1}) = \left(\frac{p_{s}}{\tau}\,, v \right) + (f_{s+1},v) \,.
\end{equation}

It is worth noting that the existence and uniqueness of the solution $p_{s+1}^{n+1}$ to the linearized equation \eqref{r1el} are proven in \cite{rh2}.  The Picard iterative method converges to a limit as $n$ approaches $\infty\,;$ and a theoretical proof of this limit is given in \cite{sdt22}.  In terms of numerical simulation, this procedure is completed when it reaches a provided halting criterion at some Picard iterative step $\alpha$th.  In order to advance to the next time step in Eq.\ \eqref{r1ed}, the previous time data is set to 
\begin{equation}\label{pdata}
	p_{s+1} = p_{s+1}^{\alpha}\,.
\end{equation}
We propose a termination indicator based on the relative successive difference throughout this work, which indicates that given a user-defined tolerance $\delta_0 > 0\,,$ the iterative process is stopped when 
\begin{equation}\label{pt}
\dfrac{\|p_{s+1}^{n+1} - p_{s+1}^{n} \|_{L^2(\Omega)}}{\| p_{s+1}^{n} \|_{L^2(\Omega)}} \leq \delta_0\,.
\end{equation}
The value $\delta_0 = 10^{-6}$ is selected in Section \ref{numer}. 

\begin{figure}[h!]
	\begin{center}
	\begin{minipage}[h]{0.27\linewidth}
			\center{\includegraphics[width=0.9\linewidth]{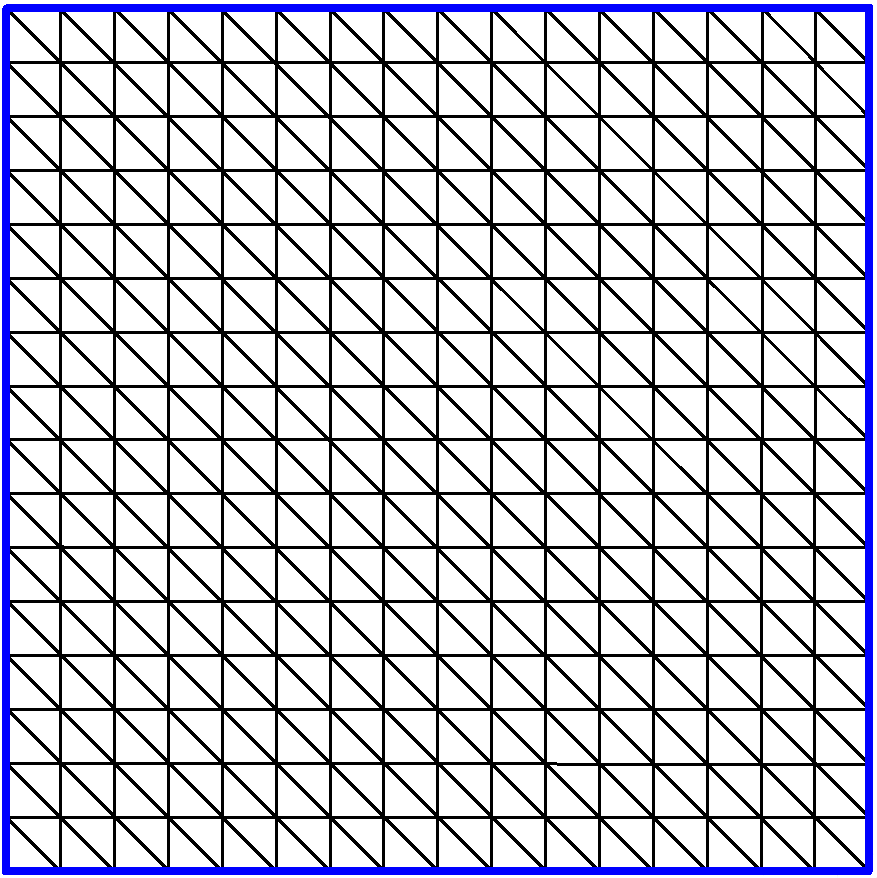}}
		\end{minipage}
	\end{center}
	\caption{Illustration of a triangular fine grid in one coarse block.}
	\label{fig:fine1block}
\end{figure}

We now discuss the fine-grid notation.  Toward discretizing \eqref{r1ed}--\eqref{r1elr}, let $\mathcal{T}_h$ (\textit{fine grid}) represent a conforming triangular partition for the computational domain $\Omega\,,$ with local grid sizes $h_P: = \tu{diam}(P)$ for all $P \in \mathcal{T}_h\,,$ and $h:= \displaystyle \max_{P \in \mathcal{T}_h} h_P\,.$  It is expected that the magnitude $h$ is sufficiently small so that the fine-grid solution is fairly close to the precise solution. The fine grid has equal squares in numerical simulation, and we use a fine-scale triangulation in which two triangles $K^h_j$ are constructed inside every square cell, as Fig. \ref{fig:fine1block} illustrates, for $j = 1,\ldots,N_c^h\,,$ with $N_c^h$ stands for the number of fine triangles in $\Omega\,.$  It should be noted that $\mathcal{T}_h$, the fine grid, is only utilized to solve local problems numerically.

Next, we denote by $V_h$ the first-order Galerkin (conventional) finite element basis space for the fine grid $\mathcal{T}_h\,:$
\begin{equation}\label{vh}
V_h:= \{ v \in V: v|_{K^h} \in \mathcal{P}_1(K^h) \; \forall K^h \in \mathcal{T}_{h}\}\,,
\end{equation}
in which $\mathcal{P}_1(K^h)$ is the space of all linear functions (that is, polynomials of degree at most 1) in the triangle $K^h\,.$ Alternatively stated, the conforming piecewise linear functions $v$
defined on the mesh $\mathcal{T}_h$ make up the finite element space $V_h\,.$ 
The initial $p_{h,0} = P_h p_0$ in this $V_h$ is given by $P_h\,,$ which is the $L^2(\Omega)$ projection operator onto $V_h\,,$ and $p_0$ comes from \eqref{ini}. 

Here is the completely Picard discrete scheme in $\mathcal{T}_h\,:$ beginning with an initial $p_{h,0} \in V_h$ and guessing $p^0_{h,s+1}\in V_h$ at the temporal step $(s+1)$th, we iterate in $V_h$ to obtain $p^{n+1}_{h,s+1}\in V_h\,:$
\begin{align}\label{r1elh}
	\left(\frac{p^{n+1}_{h,s+1}}{\tau} , v \right) + a(p^{n+1}_{h,s+1},v;p^n_{h,s+1})= \left(\frac{p_{h,s}}{\tau} , v \right) + (f_{h,s+1},v)\,,
\end{align}
with $n=0,1,2, \ldots$ and all $v \in V_h$, until attaining \eqref{pt} at some $\alpha$th Picard iteration. In order to go to the subsequent time step in \eqref{r1ed}, we utilize \eqref{pdata} to determine the previous time data $p_{h,s+1} = p_{h,s+1}^{\alpha}\,.$

\section{Coarse-scale approximation}\label{coarseapp}
The purpose of this section is to determine the coarse-grid approximation of the solution to the single-continuum nonlinear Richards equation \eqref{eq:original0} (or \eqref{r1e}).




\subsection{Overview}\label{gmsover}

To reach the coarse-scale approximation, our goal is to build the offline and especially online multiscale spaces, employing the generalized multiscale finite element method (GMsFEM \cite{G1}) for the nonlinear Richard equation \eqref{r1e}.  Motivated by \cite{gne, gnone, tfcmm}, we can construct those multiscale spaces with respect to this nonlinearity by using the linearized equation \eqref{r1elh}, in which the nonlinearity is viewed as constant at each Picard iterative step (following temporal discretization). 

\begin{figure}[h!]
	\begin{center}
		\center{\includegraphics[width=0.6\linewidth]{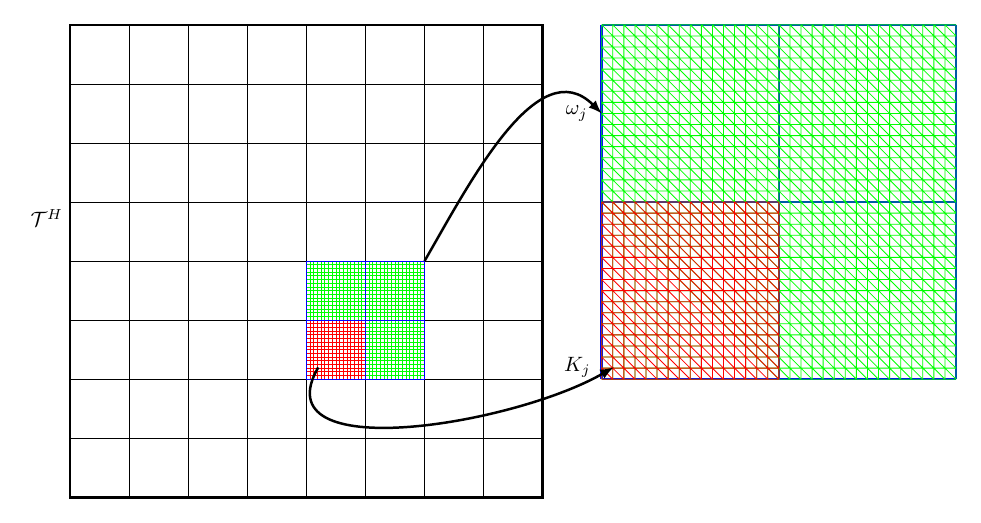}}
	\end{center}
\caption{Coarse grid $\mathcal{T}^H$ (left), coarse block $K_j$ (right, red), and coarse neighborhood $w_j$ (right).}\label{fig:coarse}
\end{figure}

First, the coarse-scale notation is covered. Consider a coarse grid $\mathcal{T}^H$ with the refinement $\mathcal{T}_h\,.$  Such coarse grid $\mathcal{T}^H$ is a typical conforming partition of $\Omega$ into equal squares (or more generally, finite elements).  Each square coarse block's edge length is represented by $H\,,$ which stands for the coarse-grid size (where $h \ll H$).  In $\mathcal{T}^H\,,$ every element is called a coarse block (or patch, or element). 
We designate $N_p$ as the total number of coarse blocks and $N_v$ as the total number of vertices of $\mathcal{T}^H$ (including those on the boundary).  The set of vertices (nodes) in $\mathcal{T}^H$ is denoted by $\{\bfa{x}_j\}^{N_v}_{j=1}\,.$ The union of all coarse blocks $K_m \in \mathcal{T}^H$ having such coarse-grid node $\bfa{x}_j$ defines the $j$th local domain or coarse neighborhood of  $\bfa{x}_j$ as follows (see Fig.~\ref{fig:coarse} as well):
\beq
\omega_j = \bigcup \{K_m \in \mathcal{T}^H : \bfa{x}_j \in \overline{K_m}\}\,.
\eeq  

We refer the readers to \cite{G1, G2, chungres1, chungres, chung2016adaptive} for the GMsFEM's details and to \cite{mcl,gne, mcontinua17} for a summing up. Using the GMsFEM, our major goal is to seek a multiscale solution $p_{\tu{ms}}$ that approximates the fine-scale solution $p_h$ from \eqref{r1elh} as closely as possible.
In order to achieve this, we utilize the GMsFEM \cite{G1} on a coarse grid at each Picard iteration to tackle local problems (which will be explained later) in every coarse neighborhood.
This enables us to consistently produce multiscale basis functions (the solution's degrees of freedom) while preserving fine-scale features.  It is evident that a key component of the GMsFEM is the construction of local multiscale basis functions.  
First, in accordance with \cite{gne}, only the so-called offline multiscale basis functions (which are calculated in the offline stage) will be generated.
Secondly, to improve the accuracy of the multiscale approximation, we will locally construct additional online multiscale basis functions (that are problem-dependent), based on such offline basis functions and several local residuals. In one of the authors' previous work \cite{denisonr}, its results show that the multiscale solution $p_{\tu{ms}}$ will quickly converge to the fine-scale solution $p_h$ when offline and online basis functions are combined in the final multiscale space $V_{\tu{ms}}\,,$ as with \cite{chungres,gne}.  For the GMsFEM, we note that our equation \eqref{eq:original0} 
possesses multiscale high-contrast coefficient $\kappa$ (depending on $\bfa{x}$).
\subsection{Offline GMsFEM}\label{offstage}
 
Now, we introduce a technique for generating offline multiscale basis functions (mostly inspired by \cite{rtt21,mcl} based on \cite{G1}).  In short, a local snapshot space is created for each coarse neighborhood $\omega_j$ using the GMsFEM's framework \cite{G1} as described above. Next, we solve an appropriate local spectral problem (specified on the snapshot space) to establish an  offline multiscale function space $V_{\tu{off}}\,.$ 

That is to say, in the following Sections \ref{snaps}, \ref{mssp}, 
provided $p_{\tu{off},s}$ (at the $s$th time step) and $p^{n}_{\tu{off},s+1}$ (at the $(s+1)$th time step  and the $n$th Picard iteration), we will build $V_{\tu{off}}= V^n_{\tu{off},s+1}\,.$  Both the snapshot functions and the multiscale basis functions in this case are independent of time.  In practice (Section \ref{offgms}), provided $p_{\tu{ms},0}$ and a starting guess $p^{0}_{\tu{ms},1}\,,$ we will need to construct only one $V_{\tu{off}} =V^0_{\tu{off},1}$ ($n=0,\ s+1 = 1$).

%


\subsubsection{Snapshot spaces}\label{snaps} 
Let $V_h(\omega_j)$ be the fine-scale approximation (FEM) space, which is the restriction of the conforming space $V_h$ to the local domain $\omega_j\,.$
The collection of all the nodes of the fine grid $\mathcal{T}_h$ situated on $\partial \omega_j$ is denoted by $J_h(\omega_j)$.  Also, $J_h(\omega_j)$ has the cardinality represented by $N_{J_j}\,.$ The index $k$ is allowed to vary as follows: $1 \leq k \leq N_{J_j}\,.$ 

We first solve the following local snapshot problem on every coarse neighborhood $\omega_j\,:$ seek the $k$th snapshot function $\phi_{k}^{j,\tu{snap}} \in V_h(\omega_j)$ such that
\beq\label{snapu}
\bsp
-\div ( \varkappa(\bfa{x},p^n_{\tu{off},s+1})\nabla \phi_{k}^{j,\tu{snap}} ) &= 0 \ \ \ \text{in} \ \omega_j,\\
\phi_{k}^{j,\tu{snap}} &= \delta_{k} \ \ \ \text{on} \ \partial \omega_j\,,
\end{split}
\eeq
where the 
 function $\delta_{k}$ is defined as
\[\delta_{k}(\bfa{x}^j_m) =
\begin{cases}
1 \quad m = k\,,\\
0 \quad m \ne k\,,
\end{cases}
\]
over all $\bfa{x}^j_m$ in $J_h(\omega_j)\,,$ for $1 \leq k \leq N_{J_j}\,.$  
Consequently, we derive the $j$th local snapshot space 
\[V_{\tu{snap}}(\omega_j) = \text{span}\{ \phi_{k}^{j,\tu{snap}} \, \bigr | \,  1 \leq k \leq N_{J_j} \}\,.\]

\subsubsection{Offline multiscale space}\label{mssp}
In order to construct multiscale basis functions on $\omega_j\,,$ local spectral problems are solved as follows: identify the $k$th eigenfunction $\psi_{k}^{j} \in V_{\tu{snap}}(\omega_j)$ and the real eigenvalue $\lambda_{k}^{j}$ that corresponds to it, such that for every $\xi$ in $V_{\tu{snap}}(\omega_j)\,,$   
\beq\label{eeunc}
a^{j}(\psi_{k}^{j},\xi) = \lambda_{k}^{j} s^{j}(\psi_{k}^{j},\xi)\,.
\eeq
These operators have the following definitions for all $\phi, \psi \in V_{\tu{snap}}(\omega_j)$ \cite{gne,mcontinua17, mcl}:  
\beq\label{schi}
\bsp
a^{j}(\phi,\psi) = \int_{\omega_j} \varkappa(\bfa{x},p^n_{\tu{off},s+1}) \nabla \phi \cdot \nabla \psi\dx,\\
s^{j}(\phi,\psi) = \int_{\omega_j} \varkappa(\bfa{x},p^n_{\tu{off},s+1})\left( \sum_{l=1}^{N_v} |\nabla \chi_{l}|^2\right) \phi \, \psi  \, \dx\,.
\end{split}
\eeq
Here, for the coarse-grid node $\bfa{x}_l\,,$ each $\chi_{l}$ represents a 
standard coarse-grid nodal basis function \cite{denisonr}, that is, with linear boundary conditions for cell problems \cite{pou}. It is noted that $\{\chi_{l}\}_{l=1}^{N_v}$ is a set of partition of unity functions for the coarse grid $\mathcal{T}^H\,.
$  
Specifically, $\forall K_m \in \omega _j\,,$ each $\chi_{j}$ is computed as follows \cite{pou}:
\begin{equation}\label{poueq}
\begin{split}
     -\div (\nabla \chi_{j}) &= 0 \ \  \mbox{in} \ \ K_m \in \omega _j \,,  \\
     \chi_{j} &= \chi^0_{j} \ \  \mbox{on} \ \ \partial K_m\,,   
\end{split}   
\end{equation}
where each $\chi^0_j$ is a standard linear partition of unity function, that is linear and continuous in $w_j\,,$ 
has value 1 at the center $\bfa{x}_j$ and 0 at the remaining coarse vertices of $w_j\,.$

%
The first $L_{\omega_j}$ eigenfunctions are obtained by sorting the eigenvalues $\lambda_{k}^{j}$ from \eqref{eeunc} in ascending order. These eigenfunctions are still designated as $\psi_{1}^{j}, \ldots, \psi_{L_{\omega_j}}^{j}\,.$   Now, for $1 \leq k\leq L_{\omega_j}$ ($\leq N_{J_j}$) and $\{\chi_{l}\}_{l=1}^{N_v}$ given in \eqref{poueq}, the $k$th local offline multiscale basis function on $\omega_j$ is defined by 
\begin{equation}\label{mbsu}
\psi_{k}^{j,\tu{off}} = \chi_{j} \psi_{k}^{j}\,.
\end{equation}

The definition of local auxiliary offline multiscale space is
\begin{equation}\label{localmbf}
	V^n_{\tu{off},s+1}(\omega_j) = \text{span}\left \{\psi_{k}^{j,\tu{off}} \, \bigr | \, 1 \leq k\leq L_{\omega_j} \right \}\,.
\end{equation}
Then, the global offline multiscale function space is obtained by 
\begin{equation}\label{globalmbf}
	V^n_{\tu{off},s+1} = \sum\limits_{j=1}^{N_v} V_{\tu{off}}(\omega_j)= \text{span}\left \{\psi_{k}^{j,\tu{off}} \, \bigr | \, 1 \leq j \leq N_v \,, 1\leq k\leq L_{\omega_j} \right \}\,,
\end{equation}
which will be used in the upcoming $(n+1)$th Picard iteration to seek the solution.

\subsubsection{Offline GMsFEM for nonlinear Richards equation}\label{offgms}
As in \cite{rtt21,gne, cemnlporo}, our strategy is to solve the problem (\ref{eq:original0}) with equivalent variational formula \eqref{r1e} by employing the offline GMsFEM (from Sections \ref{snaps}, \ref{mssp}) and the correspondingly established offline multiscale function space $V_{\tu{off}}= V^0_{\tu{off},1}$ (for $n=0$ and $s+1=1$ in \eqref{globalmbf}), during each Picard iteration of linearization at the fixed time step $(s+1)$th.  
Given an initial $p_{\tu{off},0} = P_{\tu{off}} p_0$ in $V_{\tu{off}}$ where $P_{\tu{off}}$ is the $L^2(\Omega)$ projection operator onto $V_{\tu{off}}$ and $p_0$ is from \eqref{ini}.  
In numerical simulations, for simplicity, we can take $p_0 = 0\,.$
The model reduction algorithm used through the online stage (with the source term) is as follows: given $p_{\tu{off},s}\,,$ a guess $p^0_{\tu{off},s+1} \in V_{\tu{off}}$ at the temporal step $(s+1)$th, and an obtained $p^n_{\tu{off},s+1}$ (with $n = 0, 1, 2, \ldots$), we iterate in $V_{\tu{off}}$ to solve for $p^{n+1}_{\tu{off},s+1}\, :$  
\begin{align}\label{ongmspicard}
\begin{split}
\left(\frac{p^{n+1}_{\tu{off},s+1}}{\tau} , v \right) + a(p^{n+1}_{\tu{off},s+1},v;p^n_{\tu{off},s+1}) 
= \left(\frac{p_{\tu{off},s}}{\tau} , v \right) + (f_{s+1},v)\,,
\end{split}
\end{align}
with $v \in V_{\tu{off}}$, for $n=0,1,2, \ldots\,,$ until attaining \eqref{pt} at some $\alpha$th Picard iterative step.  
Then, \eqref{pdata} is employed to select the prior time data $p_{\tu{off},s+1} = p_{\tu{off},s+1}^{\alpha}$ to progress to the subsequent temporal step in \eqref{r1ed}.


\subsection{Online GMsFEM}\label{onstage}

Given the nonlinear nature of our problem, we must consider the coefficients' changing over time. While updating offline multiscale basis functions at certain time steps might be taken into account, this method comes with a high computational cost. Including at least one online multiscale basis function in each local region (which reduces the residual) would be far more profitable. This approach follows the residual-driven online GMsFEM \cite{chungres}. The process that is being presented will enable us to quickly reduce errors based on local residual information (by adding online multiscale basis functions) while notably reducing the number of offline multiscale basis functions. That is, we will be able to obtain an improved accuracy with fewer multiscale basis functions since they incorporate the coefficient $\varkappa$'s variations over time.

Our big picture is that given $p_{\tu{off},s}$ at the previous $s$th time step, then at the current time step $(s+1)$th, we obtain the previous data $p^{n}_{\tu{off},s+1}$ from \eqref{ongmspicard} and the final multiscale solution $p^{n+1}_{\tu{off},s+1}$ (at the Picard iteration $(n+1)$th) within the obtained offline space $V_{\tu{off}}\,.$  Next, employing these $p_{\tu{off},s}\,,$ $p^{n}_{\tu{off},s+1}$ and $p^{n+1}_{\tu{off},s+1}\,,$ we compute the online multiscale basis function $\Phi^n_{j,s+1}$ over the local domain $w_j$ and get $V^n_{\tu{ms},s+1} = V_{\tu{off}} \oplus \text{ span} \{ \Phi^n_{j,s+1} \}$ (at the Picard iteration $n$th).
Last, with the previous solution $p^n_{\tu{off},s+1} \in V_{\tu{off}}\,,$ we make one more Picard iteration \eqref{r1el} to obtain the new solution $p^{n+1}_{\tu{ms},s+1} \in V^n_{\tu{ms},s+1}$ at the $(n+1)$th Picard step. 

More specifically, at this online stage, let $V_j = H_0^1(w_j) \cap V_h\,.$ 
 Using the procedure in \cite{denisonr,gne}, which is derived from the residual-based online basis enrichment for GMsFEM \cite{chungres}, our online multiscale basis function $\eta^n_{j,s+1} \in V_j$ is the solution to the following local problem in $w_j\,:$
\begin{equation}\label{ron}
d_n^j(\eta^n_{j,s+1},v) = R^j(v), \quad \forall v \in V_j\,,
\end{equation}
where
\begin{equation}\label{ain}
d_n^j(\eta^n_{j,s+1},v) = \int_{w_j} \frac{\eta^n_{j,s+1}}{\tau} v \dx + \int_{w_j} \varkappa(p^{n}_{\tu{off},s+1})\nabla \eta^n_{j,s+1} \cdot \nabla v \, \dx \,, 
\end{equation}
and
\begin{align}\label{rj}
R^j(v) = \left(\int_{w_j} \frac{p_{\tu{off},s}}{\tau}\, v \dx + \int_{w_j} fv \dx \right) - \left( \int_{w_j} \frac{p^{n+1}_{\tu{off}, s+1}}{\tau}\, v \dx + \int_{w_j} \varkappa(p^{n}_{\tu{off},s+1})\nabla p^{n+1}_{\tu{off},s+1} \cdot \nabla v \, \dx \right)\,,
\end{align}
with zero Dirichlet boundary condition $\eta^n_{j,s+1}=0$ on $\partial w_j\,.$  
We equip the space $V_j$ with the enery norm $\| \cdot \|_{V_j}\,,$ for any $v\in V_j\,:$
\begin{equation}\label{vnorm}
\|v\|^2_{V_j} = d^j_n(v,v) = \int_{w_j} \frac{v^2}{\tau}  \dx + \int_{w_j} \varkappa(p^{n}_{\tu{off},s+1}) |\nabla v|^2 \, \dx\,.
\end{equation}
The norm of the operator $R^j$ is represented by $r_j$ and computed by
\begin{equation}\label{rnorm}
r_j = \|R^j\|_{V_j} = \|\eta^n_{j,s+1}\|_{V_j}\,,
\end{equation}
introduced in \eqref{vnorm}.

Now, the online multiscale basis function over the local domain $w_j$ is of the form 
\begin{equation}\label{Phion}
\Phi^n_{j,s+1} = \chi_j \eta^n_{j,s+1}\,.
\end{equation}
Here, for every coarse vertex $\bfa{x}_j\,,$ each $\chi_{j}$ denotes a 
standard coarse grid nodal basis function, and $\{\chi_{j}\}_{j=1}^{N_v}$ is a set of partition of unity functions (for the coarse grid $\mathcal{T}^H$) defined in \eqref{poueq}. The new multiscale space is
\begin{equation}\label{Vmson}
V^n_{\tu{ms},s+1} = V_{\tu{off}} \oplus \text{ span} \{ \Phi^n_{j,s+1} \}\,.
\end{equation}

Note that at the $n$th Picard step and the $(s+1)$th time step, to add more online multiscale basis functions over each coarse neighborhood, one can iterate the online procedure using residual calculation above.  Nevertheless, in our paper, it is sufficient to add only one online multiscale basis function per coarse neighborhood.

Next, provided $p_{\tu{off},s}\,,$ using the previous solution $p^n_{\tu{off},s+1} \in V_{\tu{off}}\,,$ we make one more Picard iteration \eqref{ongmspicard} at the $(n+1)$th Picard step in $V^n_{\tu{ms},s+1}$ \eqref{Vmson} to obtain the new solution $p^{n+1}_{\tu{ms},s+1} \in V^n_{\tu{ms},s+1}$ as follows (with $v \in V^n_{\tu{ms},s+1}$):
\begin{align}\label{ongmspicard1}
\begin{split}
\left(\frac{p^{n+1}_{\tu{ms},s+1}}{\tau} , v \right) + a(p^{n+1}_{\tu{ms},s+1},v;p^n_{\tu{off},s+1}) 
= \left(\frac{p_{\tu{off},s}}{\tau} , v \right) + (f_{s+1},v)\,.
\end{split}
\end{align} 

In practice, the online multiscale basis function space $V^n_{\tu{ms},s+1}$ can be enriched (directly or by DNN prediction) not for every time step, but every 5 time steps ($s+1 = 1,5,10,15,20$ in our experiments shown in Section \ref{numer}), for example.   

Let 
\begin{equation}\label{nce}
N_c = \tu{dim}(V^n_{\tu{ms}, s+1})\,.
\end{equation}

\bigskip




\textbf{For implementation:} 

The multiscale basis functions in $V^n_{\tu{ms}, s+1}$ are represented by $\psi_l$ ($l=\overline{1,N_c}\,,$ where $N_c$ is defined in \eqref{nce}). For the sake of simplicity, we denote $p^{n+1}_{\tu{ms}, s+1}$ by $p^{n+1}_{c,s+1}\,,$ 
and the subscript $c$ in \eqref{matrixms} 
is dropped in the following.  The approximate solution $p^{n+1}_{c,s+1}$ at the time step $(s+1)$th is 
\begin{equation}\label{pcapprox}
p^{n+1}_{c,s+1} \approx \sum_{l=1}^{N_c} p^{n+1}_{c,s+1}(\bfa{x}_l) \, \psi_l = \sum_{l=1}^{N_c} p^{n+1}_{c,s+1,l} \, \psi_l\,,
\end{equation}
in which $p^{n+1}_{c,s+1,l}$ is the point value of $p^{n+1}_{c,s+1}$ at $\bfa{x}_l$ with $l \in N_c\,.$  
Let us denote 
\begin{equation}\label{pcoord}
\bfa{p}^{n+1}_{c,s+1} = (p^{n+1}_{c,s+1,1}, \ldots, p^{n+1}_{c,s+1,N_c})\,.
\end{equation}

Using the form \eqref{pcapprox} and \cite{richarde1}, we can write
\begin{align}\label{comktc}
\begin{split}
\varkappa(\bfa{x},p^n_{\tu{off},s+1}) & \approx \sum_{l=1}^{N_c} \varkappa(\bfa{x},p^n_{\tu{off},s+1,l})\psi_l = \sum_{l=1}^{N_c} \varkappa_l \, \psi_l\,,\\
p_{\tu{off},s}  & \approx \sum_{l=1}^{N_c} p_{\tu{off},s} (\bfa{x}_l) \psi_l = \sum_{l=1}^{N_c} p_{\tu{off},s,l} \, \psi_l\,, \\
p^n_{\tu{off},s+1} &\approx \sum_{l=1}^{N_c} p^n_{\tu{off},s+1}(\bfa{x}_l) \psi_l = \sum_{l=1}^{N_c} p^n_{\tu{off},s+1,l} \, \psi_l\,.
\end{split}
\end{align}


For $l,r = \overline{1,N_c}\,,$ we let
\begin{align}\label{comstiffc} 
\bfa{A}^n_{c,s+1} = \{A^n_{c,s+1,lr}\}\,, \quad
\bfa{b}_{c,s+1} = \{b_{c,s+1,l}\}\,, 
\end{align}
with
\begin{align}\label{comstiffc2}
\begin{split}
A^n_{c,s+1,lr} &= \int_{\Omega} \varkappa(\bfa{x},p^n_{\tu{off},s+1}) \nabla \psi_l \cdot \nabla \psi_r \, \dx \,,\\
b_{c,s+1,l} &= \int_{\Omega} \left(\frac{p_{\tu{off},s}}{\tau} + f_{s+1}\right)\, \psi_l \, \dx \,.
\end{split} 
\end{align}

\bigskip

Employing the subscript $c$ again in \eqref{comstiffc}, we can rewrite Eq.\ \eqref{ongmspicard1} in the following matrix form:
\begin{equation}\label{matrixms} \frac{\bfa{p}^{n+1}_{c,s+1}}{\tau} + \bfa{A}^n_{c,s+1} \ \bfa{p}^{n+1}_{c,s+1} = \bfa{b}_{c,s+1}\,,
\end{equation} 
where $\bfa{A}^n_{c,s+1}$ is the coarse-scale stiffness matrix, $\bfa{b}_{c,s+1}$ is the load vector, and $\bfa{p}^{n+1}_{c,s+1}$ is the vector of $N_c$ components.  
 
Note that our solution $p^{n+1}_{\tu{ms},s+1}$ is in the continuous space $V^n_{\tu{ms},s+1}\,.$  Thus, in numerical computation, at the current Picard step $n$th, we employ the first-order finite elements on the fine mesh $\mathcal{T}_h$ (through the standard FEM basis) to express the multiscale basis functions.  Also, we denote by $\bfa{A}^n_{s+1}$ and $\bfa{b}_{s+1}$ the fine-scale stiffness matrix and right-hand side vector computed by the finite element basis functions at the nodes of the fine grid $\mathcal{T}_h\,.$ 
If all the multiscale basis functions are collected and their fine-scale coordinate representation is arranged in columns, then these columns form the downscaling operator $\bfa{R}_g\,.$ 
Hence, instead of using \eqref{comstiffc2}, we can utilize the fine-scale forms $\bfa{A}^n_{s+1}\,, \bfa{b}_{s+1}$ 
to compute
\begin{align}\label{calcoarse}
\bfa{A}^n_{c,s+1} = \bfa{R}_g^{\tu{T}} \bfa{A}^n_{s+1} \bfa{R}_g\,, \qquad \bfa{b}_{c,s+1} = \bfa{R}_g^{\tu{T}} \bfa{b}_{s+1}\,.
\end{align}

Finally, the coarse-grid multiscale solution $p^{n+1}_{\tu{ms},s+1}$ possesses the fine-grid coordinate representation $\bfa{p}^{n+1}_{\tu{ms},s+1}$ in the form of vector
\begin{equation}\label{finerepsol}
\bfa{p}^{n+1}_{\tu{ms},s+1} = \bfa{R}_g \bfa{p}^{n+1}_{c,s+1}\,.
\end{equation}


\section{Deep learning for online GMsFEM}\label{dlgms}



\begin{figure}[h!]
  \centering
\begin{tikzpicture}
\centering
\node (input) [Input] {\textbf{Given input data} \\ \begin{itemize}[leftmargin=.2in] \item Generate permeability field realizations
  \item Compute online multiscale basis functions
\end{itemize}};
\node (dnn) [DNN, right = 1cm of input] {\textbf{Deep neural networks} \noindent \begin{itemize}[leftmargin=.2in] \item Construct and train DNNs \item Predict new online multiscale basis functions \end{itemize}};
\node (sm) [MS, right = 1cm of dnn] {\textbf{Multiscale method} \\
\begin{itemize}[leftmargin=.2in] \item Solve for coarse-scale solution \end{itemize}} ;
\draw [arrow] (input) -- (dnn);
\draw [arrow] (dnn) -- (sm);
\end{tikzpicture}
\caption{A flowchart depicting the concept of utilizing deep learning in the online GMsFEM framework.}
\label{fig:flowcharts}
\end{figure}
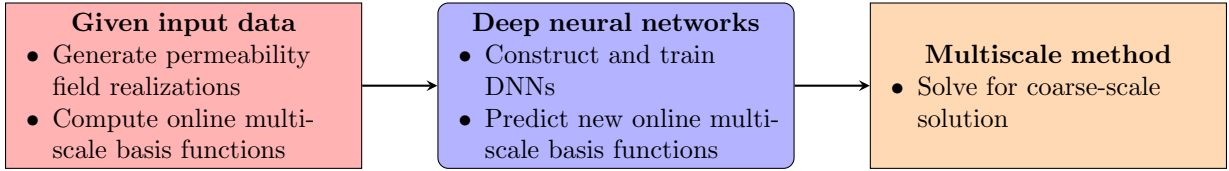

Within applications, the heterogeneous permeability field $\kappa(\bfa{x})$ in the flow problem \eqref{eq:original0} has uncertainties in specific local regions. To calculate the uncertainties of the flow solution, thousands of forward simulations are required. For computing the solutions efficiently and precisely, the online GMsFEM offers us a rapid solver.
Given a large volume of simulation data, it is interesting to develop a technique which makes use of the provided data, in order to minimize the amount of direct computational work that needs to be done then. 
For the Richards equation \eqref{r1e}, it is assumed that the Picard linearization procedure \eqref{r1elr} ceases at $(n+1)$th Picard iteration and the $(s+1)$th temporal step.
Our goal in this study is to use deep neural networks (DNNs) to simulate the connection between the online GMsFEM's main component, namely online multiscale basis functions, and the heterogeneous permeability coefficients $\kappa(\bfa{x})$ (so the hydraulic conductivity coefficients $\varkappa(\bfa{x},p(t,\bfa{x}))$).  It should be noted that the Karhunen-Lo\`eve expansion (KLE) \cite{sdt22,cnn-mcmc20} is used to create our random heterogeneous permeability fields. 
 After this relationship is established, we can supply the network with any realization of the permeability field (so the conductivity field with the provided solution data from the previous Picard iteration) to extract the matching online multiscale basis functions and then restore the fine-scale GMsFEM solution to \eqref{eq:original0}. In the remainder of this work, for the sake of simplicity, we can omit the subscript $(s + 1)$ when there are Picard steps. Fig.~\ref{fig:flowcharts} provides an overview of the general concept of using deep learning in the online GMsFEM framework.

Our idea is that at the time step $(s+1)$th and Picard iteration $n$th, over each coarse neigborhood $w_j\,,$ correspondingly to the local online multiscale basis function $\Phi^{n}_{j,s+1}$ \eqref{Phion}, we will predict a local online multiscale basis function $\Phi^{n,\tu{pred}}_{j,s+1}$ by the following neural network and obtain the space 
\begin{equation}\label{vpred}
V^{n,\tu{pred}}_{\tu{ms},s+1} = V_{\tu{off}} \, \oplus \, \tu{span} \{\Phi^{n,\tu{pred}}_{j,s+1}\}\,.
\end{equation}
Subsequently, we compute the predicted solution $p^{n+1,\tu{pred}}_{\tu{ms},s+1} \in V^{n,\tu{pred}}_{\tu{ms},s+1}$ from \eqref{finesolpred} and compare this prediction $p^{n+1,\tu{pred}}_{\tu{ms},s+1}$ with the solution $p^{n+1}_{\tu{ms},s+1} \in V^n_{\tu{ms},s+1}$ obtained from \eqref{ongmspicard1}. 

Given a permeability field $\kappa(\bfa{x})$ and the solution $p^n_{\tu{off},s+1}$ at the previous $n$th Picard step and the $(s+1)$th time step, thus given a hydraulic conductivity field \[\varkappa(\bfa{x},p^n_{\tu{off},s+1}) = \kappa(\bfa{x})\, \mu(p^n_{\tu{off},s+1}) = \kappa(\bfa{x}) \cdot \frac{1}{1+|p^n_{\tu{off},s+1}|}\,,\]
 which is abbreviated by $\varkappa(p^n_{\tu{off},s+1})\,,$
 one can compute the multiscale basis functions $\Phi^n_{j,s+1}$ defined by \eqref{Phion} for $1 \leq j \leq N_v\,.$ 

Theoretically, it is intriguing to create the nonlinear map $g_{B}^{j,n}\,,$ which maps the permeability coefficient $\kappa$ to the expected local online multiscale basis function $\Phi^n_{j,s+1}(\kappa)\,,$ over the $j$th coarse neighborhood $\omega_j$ (for $1 \leq j\leq N_v$):
\begin{equation}\label{gb}
g_{B}^{j,n}: \kappa \to \Phi^n_{j,s+1}(\kappa)\,.
\end{equation}
This map serves only as inspiration.  More significantly, in practice, our nonlinear neural network below will numerically approximate it. 

Our research aims to solve the heterogeneous nonlinear Richards equation \eqref{eq:original0} as fast and precisely as possible, by using deep learning to rapidly approximate the above desired online multiscale basis function related to the uncertainties of the permeability field $\kappa\,.$  With every given realization $\kappa\,,$ its image $\Phi^n_{j,s+1}(\kappa)$ under the local online multiscale basis map $g_{B}^{j,n}$ can be computed by \eqref{Phion} in Section \ref{onstage}.  We use these forward calculations as training data to construct a deep neural network that will roughly represent the corresponding map, that is, 
\begin{align}\label{nnbm}
\mathcal{N}^{j,n}_{B}(\kappa) &\approx g^{j,n}_{B} (\kappa)\,.
\end{align}  

Such neural network is designed to predict local online multiscale basis function $\Phi^n_{j,s+1}$ as the output, with the permeability field $\kappa$ as the input. After establishing this neural network, over the related coarse neighborhood, local online multiscale basis function can be predicted for every new permeability realization $\kappa\,.$   By utilizing this local information obtained from the neural network and the pre-computed global information, we can build the downscale operator $\bfa{R}_{\tu{pred}}\,,$ which includes all multiscale basis functions (with their fine-scale coordinate representation displayed in columns). Finally, we can solve the linear system \eqref{apred} and derive the multiscale solution through \eqref{finesolpred}.

In practice, as Section \ref{numer}, over each coarse neighborhood $w_j\,,$ only one online multiscale basis function is added.  Hence, only one neural network is constructed over $\Omega$ (that is, over all $N_v$ coarse neighborhoods $w_j$ where $1 \leq j \leq N_v$).


\subsection{Network architecture}\label{nnarch}
Deep neural networks typically consist of three groups: 
the input layer, one or a number of hidden layers, and the output layer.  These layers are arranged in a chain pattern in most deep neural network architectures, with each layer acting as a function of the previous layer \cite{rpdnn, dnngood}.
Hence, there is at least one hidden layer in any neural network (also known as artificial neural network).
Note that in research papers, each DNN's layer consists of neurons (units); however, in simulations, the input layer has data as input vector (not neurons).
Deep neural networks are defined as possessing a minimum of two hidden layers, or often a lot more.

A normal representation of a fully connected, $L$-layer feedforward deep neural network $\mathcal{N}^L$ is as follows \cite{dnngms19,rpdnn, dnngood,sdt22}:
\begin{equation}\label{lnn}
\mathcal{N}(z;\theta) = \mathcal{N}^L(z;\theta) = \sigma^L(\bfa{W}^{L} \sigma^{L-1}(\cdots \sigma^2(\bfa{W}^2 \sigma^1(\bfa{W}^1 z + \bfa{c}^1) + \bfa{c}^2) \cdots) + \bfa{c}^{L})\,,
\end{equation}
in which $z$ is the input layer, $L$ is the depth of the network expressed as the number of layers (not including the input layer), $\theta:=(\bfa{W}^1,\bfa{W}^2, \ldots, \bfa{W}^{L}, \bfa{c}^1, \bfa{c}^2, \ldots, \bfa{c}^{L})$ having the  bias vectors $\bfa{c}^l$ together with the weight matrices $\bfa{W}^l$ (where $l=1,\ldots,L$), and $\sigma^l$ are the activation functions (which are usually nonlinear as SELU or ReLU at the hidden layers). 

More precisely, this neural network's input layer is
\[\mathcal{N}^0(z) = z\,.\]
With $1 \leq l < L\,,$ the network's $l$th hidden layer is
\[\mathcal{N}^l(z) = \sigma^l (\bfa{W}^{l}\, \mathcal{N}^{l-1}(z)  + \bfa{c}^{l})\,.\]
The output layer (with linear, or identity, or no activation function) is
\[\mathcal{N}^L(z) = \bfa{W}^{L}\, \mathcal{N}^{L-1}(z)  + \bfa{c}^{L}\,.\]
%
%

\begin{figure}[h!]
\centering
\includegraphics[width=0.65\linewidth]{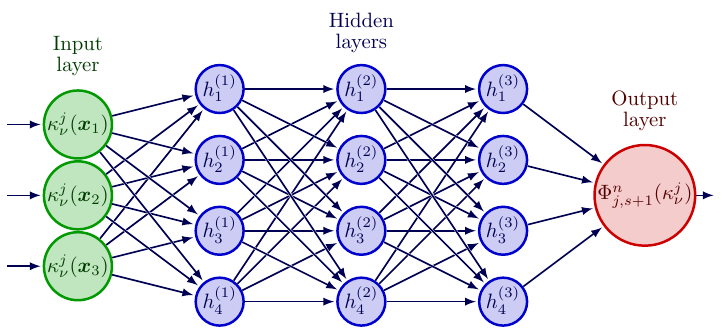}
	\caption{Visualization of a deep neural network.}
	\label{dnn_illustrate}
\end{figure}

This kind of neural network is used to approximate the desired output $y$ (which is sometimes referred to as the target,  correct output, expected output, or training value).  Our plan here is to determine $\theta^*$ through solving the optimization problem 
\begin{equation}\label{opt}
\theta^{*} = \argmin_{\theta} \mathcal{L}(\theta)\,.
\end{equation}
In this case, $\mathcal{L}(\theta)$ denotes the loss function and measures the difference between the image (predicted output) of the network $\mathcal{N}(z;\theta)$'s input and the target $y$ in a training set made up of training pairs $(z^j_{\nu},y^j_{\nu})$ (also known as training patterns, training examples, or training samples).  Here, the mean squared error serves as the loss function (under the Euclidean norm $| \cdot |$ for vectors):
\begin{equation}\label{loss}
\mathcal{L}(\theta) = \frac{1}{N_v \cdot  N_s} \sum_{j=1}^{N_v} \sum_{\nu =1}^{N_s} | y^j_{\nu} - \mathcal{N}(z^j_{\nu}; \theta)|^2\,,
%
%
\end{equation}
where $N=N_v \cdot  N_s$
represents the number of training samples for the neural network ($N_v$ is the number of coarse neighborhoods and $N_s$ is the number of training samples per each coarse neighborhood). A deep neural network example is depicted in Fig.~\ref{dnn_illustrate} (to be discussed in the following).  

Assume that we are provided a set of different realizations of permeabilities $\{\kappa^j_{\nu}\}$ with $j=1, \ldots,N_v\,, \nu=1,\ldots,N_s\,,$ for the entire domain $\Omega\,,$ where $\kappa^j_{\nu}$ is in the local domain $w_j\,.$  
Each input of our neural network is $z^j_{\nu} = \kappa^j_{\nu} \in \mathbb{R}^m\,,$ which is a vector (also known as feature vector or input vector) that contains the values of the permeability field as features $\kappa^j_{\nu}(\bfa{x}_1), \kappa^j_{\nu}(\bfa{x}_2),\ldots,\kappa^j_{\nu}(\bfa{x}_m)\,,$ where the local spatial points are $\bfa{x}_1, \bfa{x}_2,\ldots, \bfa{x}_m \in \mathbb{R}^d$ ($m$ fine-grid nodes on $w_j$), with $m=1089$ and $d=2$ in our experiments. 
It will be discussed further in Section \ref{numer} that our random permeability fields are generated by the Karhunen-Lo\`eve expansion (KLE) \cite{cnn-mcmc20,sdt22}.
Corresponding to the input $z^j_{\nu} = \kappa^j_{\nu}\,,$ the precise output $y^j_{\nu}$ is coefficient vector of the vectorized online multiscale basis function $\Phi^n_{j,s+1}(\kappa^j_{\nu})$ \eqref{Phion}, with $m$ entries of fine-scale coordinates, on the coarse neighborhood $\omega_j\,.$  By minimizing the loss function \eqref{loss} with regard to the network parameter $\theta_B\,,$ the deep neural network $\mathcal{N}_{B}^{j,n}(z, \theta_B)$ is trained using the sample pairs $(z^j_{\nu}, y^j_{\nu})\,.$ This allows the trained neural network $\mathcal{N}_{B}^{j,n}(z, \theta^*_B )$ to approximate the map $g_{B}^{j,n}$ \eqref{gb} on coarse mesh. Posteriorly to establishing the neural network, we apply it to any new permeability field $\kappa^j_{N+\iota}$ ($\iota = 1,\ldots,M\,, j=1, \ldots,N_v$) in order to rapidly predict the coarse-grid output, that is the local online multiscale basis function $\Phi^n_{j,s+1}$ \eqref{Phion} through
\begin{equation}\label{bpred}
\Phi^{n,\tu{pred}}_{j,s+1}(\kappa^j_{N+\iota}) = \mathcal{N}_{B}^{j,n}(\kappa^j_{N+\iota};\theta^*_B) \approx g_{B}^{j,n}(\kappa^j_{N+\iota}) =  \Phi^n_{j,s+1}(\kappa^j_{N+\iota})\,.
\end{equation} 


\subsection{Network-based multiscale solver}\label{nnsolve}

After the neural network is designed, benefiting from \eqref{matrixms}, at the final Picard step $(n+1)$th and last time step $(s+1)$th, we can obtain the predicted global coarse-grid solution $p^{n+1,\tu{pred}}_{c,s+1}$ of the Richards equation \eqref{eq:original0} through solving for the predicted coarse-grid coefficient vector $\bfa{p}^{n+1,\tu{pred}}_{c,s+1}$ \eqref{pcoord}
from the following linear algebraic system
\begin{equation}\label{apred} \frac{\bfa{p}^{n+1,\tu{pred}}_{c,s+1}}{\tau} + \bfa{A}^n_{c,s+1} \ \bfa{p}^{n+1,\tu{pred}}_{c,s+1} = \bfa{b}_{c,s+1}\,.
\end{equation}
Here, $\bfa{A}^n_{c,s+1}$ and $\bfa{b}_{c,s+1}$ were computed in \eqref{comstiffc2}.  Also, note that in this situation, the solution $\bfa{p}^{n+1,\tu{pred}}_{c,s+1}$ is the same as solution $\bfa{p}^{n+1}_{c,s+1}$ from \eqref{matrixms}, but their fine-scale representations are different, as follows. 

Now, our numerical implementation computes the multiscale basis functions using first-order finite elements over the fine grid $\mathcal{T}_h$ (via the FEM basis). Thus, our downscaling operator $\bfa{R}_{\tu{pred}}$ (as a matrix) is directly formed by all offline multiscale basis functions \eqref{mbsu} and the predicted coarse-grid online multiscale basis functions \eqref{bpred} (having their fine-grid coordinate representation in column vectors), that is, 
\begin{equation}\label{rpredf}
\bfa{R}_{\tu{pred}} = \left [\psi_{1}^{1,\tu{off}}, \ldots , \psi_{L_{\omega_1}}^{1,\tu{off}},
\ldots,
\psi_{1}^{N_v,\tu{off}}, \ldots , \psi_{L_{\omega_{N_v}}}^{N_v,\tu{off}},
\Phi^{n,\tu{pred}}_{1,s+1}(\kappa^1_{N+\iota}), \ldots, \Phi^{n,\tu{pred}}_{N_v,s+1}(\kappa^{N_v}_{N+\iota})\right]\,,
\end{equation}
For each coarse neighborhood, note that the number of testing samples  is $M$ ($\iota=1, \ldots, M$), so there are $M$ different operators $\bfa{R}_{\tu{pred}}\,.$  As \eqref{finerepsol}, we then obtain the predicted coarse-grid multiscale solution  $p^{n+1,\tu{pred}}_{\tu{ms},s+1}$ having its fine-grid coordinate representation $\bfa{p}^{n+1,\tu{pred}}_{\tu{ms},s+1}$ as vector
\begin{equation}\label{finesolpred}
\bfa{p}^{n+1,\tu{pred}}_{\tu{ms},s+1} = \bfa{R}_{\tu{pred}}\ \bfa{p}^{n+1,\tu{pred}}_{c,s+1}\,.
\end{equation}

\section{Numerical examples}\label{numer}

\bigskip

We consider in simulations a heterogeneous computational domain $\Omega = (0,1) \times (0,1)$ with uncertainties.  Our square coarse grid $\mathcal{T}^H$ is of the size $8 \times 8$ as Fig.~\ref{fig:coarse} (with the edge size of $H=1/8$), and its refinement fine grid $\mathcal{T}_h$ is of $128 \times 128$ squares 
with two triangles (having the size $h=\sqrt{2}/128$) per square fine cell as Fig.~\ref{fig:fine1block}.  The Picard iteration process starts with a pressure guess of $0\,.$ Except some situations (the multiscale basis functions, the norms \eqref{bochner} in Bochner space, and solutions at selected temporal steps), our solution calculation is at the last time step $S=20$ so that $S\tau = T=5 \cdot 10^{-5}$ with time step size $\tau = T/S=25 \cdot 10^{-7}$ (for the time-dependent case \eqref{tdc}).  We suppose the final Picard step is $(n+1)$th, and the Picard iterative stopping criterion over all numerical experiments is $\delta_0 = 10^{-6}\,,$ which guarantees the convergence of this linearization technique. For both the steady-state case and the time-dependent case, there can be a maximum of four Picard iterations. 
\begin{figure}[H]
	\begin{center}
	  \includegraphics[width=0.92\linewidth]{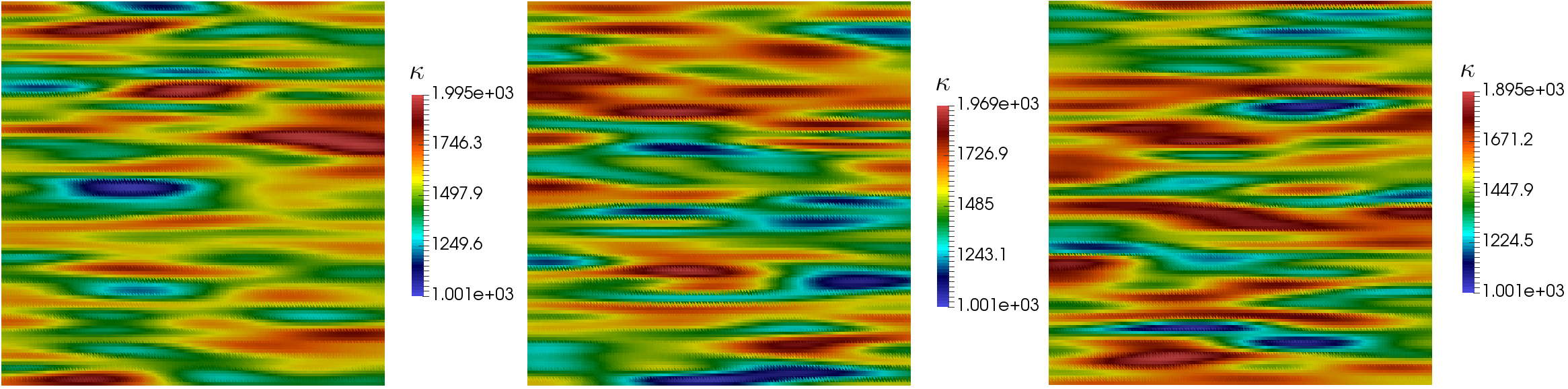}
	\end{center}
	\caption{Examples of three permeability fields $\kappa$ with values in range [10, 2000].}
	\label{fig:kappa}
\end{figure}


\subsection{Stochastic permeability}\label{sper}
As in \cite{sdt22}, we use the Karhunen-Lo\`eve expansion (KLE) \cite{aarnes2008mixed,cnn-mcmc20,cnn-homo-poro20,ganis2008stochastic,kle2} to parameterize the uncertainty of the permeability fields. Let $\Upsilon (\bfa{x},\omega)$ represent the random process.  Also, $C(\bfa{x},\hat{\bfa{x}})$ denotes its covariance function, having the position vectors $\bfa{x}=(x_1,x_2)\,, \hat{\bfa{x}}=(\hat{x}_1,\hat{x}_2)$ over $\Omega\,.$
We select a unique orthonormal basis $\{\varphi_k(\bfa{x})\}$ in the Hilbert space $L^2(\Omega)\,,$ 
where $\varphi_k(\bfa{x})$ are the eigenfunctions of $C(\bfa{x},\hat{\bfa{x}})$  in the following spectral problem:
\begin{equation}\label{spec}
\int_{\Omega} C(\bfa{x},\hat{\bfa{x}})\varphi_k(\hat{\bfa{x}})\tu{d} \hat{x} = \gamma_k \varphi_k(\bfa{x})\,, \quad k = 1,2,\ldots\,,
\end{equation}
with eigenvalues $\gamma_k\,.$  Using $\gamma_k$ and $\varphi_k$ resulted from \eqref{spec}, we may represent $\Upsilon (\bfa{x},\omega)$ in the Karhunen-Lo\`eve expansion \cite{kle2}:  
\begin{equation}\label{randfield}
    \Upsilon (\bfa{x},\omega) = \sum^{\infty}_{k=1} \sqrt{\gamma_k} \zeta_k(\omega) \varphi_k({\bfa{x})}\,,
\end{equation}
where $\zeta_k(\omega)$ are the random variables to be determined. 

To conserve most of the energy of $\Upsilon (\bfa{x},\omega)\,,$ in experiments, with respect to \eqref{randfield}, only $N_{\Upsilon}$ primary terms (attained by $\gamma_k$'s magnitude) must be maintained, hence 
\[\Upsilon (\bfa{x},\omega) = \sum^{N_{\Upsilon}}_{k=1} \sqrt{\gamma _k} \zeta_k(\omega) \varphi_k({\bfa{x})}\,.
\]
The covariance function $C(\bfa{x},\hat{\bfa{x}})$ is expected to have the following form:
\begin{equation}\label{co-matrix}
C(\bfa{x},\hat{\bfa{x}})=\sigma^2 \text{exp} \Bigg(-\sqrt{\frac{|x_1-\hat{x}_1|^2}{\eta^2_1}+\frac{|x_2-\hat{x}_2|^2}{\eta^2_2}}\Bigg),
\end{equation}represented by $\kappa(\bfa{x})\,,$ with $\eta_1 = 0.05\,, \eta_2=0.2$ as correlation lengths in each direction of space, and $\sigma^2=2$ as the variance.  Then, we can express each of our stochastic permeability fields $\kappa(\bfa{x},\omega)$ by 
\begin{equation}\label{kleform}
\kappa(\bfa{x},\omega) = \text{exp}(a_k (m(\bfa{x},\omega))),
\end{equation}
having $a_k > 0\,,$ and $m(\bfa{x},\omega ) = m(\Upsilon (\bfa{x},\omega))$ as the heterogeneous porosity. 
Recall that Fig.~\ref{fig:kappa} depicts some examples of the permeabilities we created.

It is possible to construct a stochastic permeability field in other ways without compromising our considered deep learning process.  

\subsection{Overview of numerical deep learning for the online GMsFEM}\label{numdnn}

Assume that our permeability fields are heterogeneous in $\Omega = (0,1)^2$ and contain uncertainties. 
More precisely,
we take into account a number of random realizations of distinct permeability fields $\kappa^j_1\,,\kappa^j_2\,, \ldots \,,\kappa^j_{N_s+M}$ over each coarse neighborhood $w_j\,,$ for all $j=1,\ldots,N_v\,,$ so for all $w_j$ that form $\Omega\,.$ 

With the technique and terminology provided in Section \ref{dlgms}, we create the sample pairs using the online GMsFEM in Section \ref{onstage}.  More specifically, at the $n$th Picard iteration (and at the time step $(s+1)$th for the time-dependent case \eqref{tdc}), local problems are solved to obtain the local online multiscale basis functions $\Phi^n_{j}$ \eqref{Phion} over the coarse neighborhood $w_j\,.$  Our neural network has the input as stochastic permeability field $z = \kappa$ and the expected output as local online multiscale basis function $y=\Phi^n_{j}(\kappa)$ \eqref{Phion}. At random, these sample pairs are separated into the training set (also referred to as the learning set, training data, or training dataset) and testing set.  
A large number $(N_v \cdot N_s)$ of stochastic permeability realizations $\kappa^j_1, \kappa^j_2, \ldots, \kappa^j_{N_s}$ (with $N_s$ permeability fields for each local domain $w_j$ and $j=1,\ldots,N_v$) generate sample pairs that make up the training set over all $\Omega\,.$ In contrast, the testing set, which evaluates the prediction power of the trained network, is formed by sample pairs created by the remaining $(N_v \cdot M)$ stochastic permeability realizations $\kappa^j_{N+1},\kappa^j_{N+2}, \ldots, \kappa^j_{N_s+M}$ (with $M$ permeability fields per each neighborhood $w_j$).

In \cite{dnngms19}, for each local multiscale basis function, a distinct neural network is established.  Now, in our paper, over every coarse neighborhood, we add a fix number $N_{\tu{on}}$ of online multiscale basis functions to any given set of offline multiscale basis functions.  Thus, over the whole domain $\Omega$ (consisting of all $N_v$ coarse neighborhoods $w_j$), and with all local online multiscale basis functions, we build only one 
neural network through solving an
optimization problem via minimizing the loss function defined by 
the training set. We sum up the network architecture for training local online multiscale basis functions 
as follows, with $N_v=81\,, j=1,\ldots,N_v\,,$ and $N_s=5000\,, \nu=1,\ldots,N_s\,, $ $m = 1089 = 33 \times 33\,,$ and $\bfa{x}_1,\bfa{x}_2,\ldots, \bfa{x}_m \in \mathbb{R}^2$ ($m$ fine-grid nodes on each coarse neighborhood $w_j$):
\begin{itemize}
\item For the local online multiscale basis function $\Phi^n_{j}(\kappa^j_{\nu})\,,$ we establish a network $\mathcal{N}_{B}^{j,n}$ using
\begin{itemize}
\item Input: vectorized permeability field $\kappa^j_{\nu}$ comprising values $\kappa^j_{\nu}(\bfa{x}_1), \kappa^j_{\nu}(\bfa{x}_2), \ldots, \kappa^j_{\nu}(\bfa{x}_m)\,,$
\item Desired output: coefficient vector of the local online multiscale basis function $\Phi^n_{j}(\kappa^j_{\nu})$  (with $m=1089$ entries of fine-scale coordinates with respect to $\mathcal{T}_h$) over the coarse neighborhood $\omega_j\,,$   
\item Loss function \eqref{loss}:
		mean squared error (with the Euclidean norm $| \cdot |$ on fine grid) 
  \begin{equation}\label{loss2}
  \frac{1}{N_v \cdot N_s} \sum_{j=1}^{N_v} \sum_{\nu =1}^{N_s} |\Phi_{j}^{n}(\kappa^j_{\nu}) - \mathcal{N}_{B}^{j,n}(\kappa^j_{\nu}; \theta_B)|^2\,,
  \end{equation}
\item Activation function: the first hidden layer uses SELU (Scaled Exponential Linear Unit) activation function, then all the other hidden layers employ ReLU (Rectified Linear Unit)  activation, the last output layer has the linear activation function (also known as the identity function (multiplied by 1.0), or no activation function). 

\item DNN structure: 3 hidden layers, each layer consists of 1345-1600 neurons,

\item 
Kernel initializer: ``LecunNormal'' at the first hidden layer, ``He normal'' at the second hidden layer, and ``normal'' at all the other layers (except the input layer), 

\item Training optimizer: Adam with learning rate = 0.0001.
\end{itemize}
\end{itemize}

In our experiments, over each coarse neighborhood, we only add one online multiscale basis function ($N_{\tu{on}}=1$).  Hence, only one neural network is built over $\Omega$ (that is, over all $N_v=81$ coarse neighborhoods $w_j$ where $j=1,2, \ldots,81$) for each of the numbers $Nb$ of offline multiscale basis functions per coarse neighborhood.  We choose $Nb=6\,,$ so select six cases of 2, 4, 6, 8, 12, 16 offline multiscale basis functions.  
 Therefore, six neural networks will be established in this paper.  For each neural network, the input is a matrix of size $(5000 \cdot 81) \times 1089\,,$ with $(5000 \cdot 81)$ rows of vectorized permeability fields (each of 81 coarse neighborhoods has 5000 permeability fields) and 1089 columns of fine-scale coordinates for each permeability field.  Within each neural network, the output is a matrix of size $(5000 \cdot 81) \times 1089\,,$ with $(5000 \cdot 81)$ rows of coefficient vectors of online multiscale basis functions (5000 basis functions per each of 81 coarse neighborhoods) and 1089 columns of fine-scale coordinates for each online multiscale basis function.  
 
We employ the activation functions ReLU (Rectified Linear Unit) \cite{ReLU28} and SELU (Scaled Exponential Linear Unit) \cite{selu} due to their shown efficacy for training deep neural networks in the absence of vanishing gradient issues. 
Besides, among all the nonlinear activation functions, ReLU's derivative is the simplest. 
Although SELU is less popular than ReLU, it is far more promising because it cannot die, unlike ReLU.  Self-normalization neural network (SNN) is inducible by SELU as definition, and a zero mean as well as a unit standard deviation are automatically attained by the neuronal activations in the SNN. 
%
Being widely used in neural network training, the optimizer Adam (Adaptive Moment Estimation) is a further development of stochastic gradient descent (SGD) \cite{ada29}.
%
%
We train our neural network employing Python API Tensorflow and Keras \cite{keras30} across all the experiments.

A neural network is built upon the training process, and with a given input, we can utilize it to predict the output.  The prediction accuracy plays a critical role in determining the effectiveness of this network.  In our tests, we use $(N_v \cdot M)$ sample pairs (which were not taken into account during network training) to examine the prediction ability of our created network \cite{keras30}.  That is, comparing the prediction and the precise output with respect to these sample pairs as the testing set, we compute the difference in a suitable metric, which can be summed up as follows.

Across the training set on coarse grid, the relative $l^2$ error, that is, the root mean square error (RMSE) \cite{nlnlmc31,rl2e1, rl2e2, rl2e3} is as stated in \cite{sdt22}:
%
\begin{equation}\label{rl2e}
RMSE = \sqrt{\frac{\sum_{j=1}^{N_v} \sum_{\nu =1}^{N_s}| \hat{y}^j_{\nu} - y^j_{\nu}|^2}{\sum_{j=1}^{N_v} \sum_{\nu = 1}^{N_s} |y^j_{\nu}|^2}}\,.
\end{equation}
Here, $j=1,\ldots, N_v\,,$ in which $N_v$ is the number of coarse neighborhoods $w_j$ over the whole domain $\Omega\,,$ and $N_s$ is the number of training pairs per each $w_j\,.$ We choose $N_v=81$ and $N_s=5000\,.$ Also,
$y^j_{\nu}$ is the correct output, and $\hat{y}^j_{\nu}$ is the predicted output of the training pairs $(z^j_{\nu}, y^j_{\nu})$ defined right after \eqref{opt}.   
%
%
%

With each of $M$ testing samples, at the $n$th Picard step (and at the current time step $(s+1)$th for the time-dependent case \eqref{tdc}), over the neighborhood $w_j\,,$ the relative $l^2$ error between the desired output $\Phi^n_{j}$ \eqref{Phion} (vectorized) and their corresponding predicted local online multiscale basis functions $\Phi^{n,\tu{pred}}_{j}$ (vectorized) defined in \eqref{bpred} is
\begin{equation}\label{Cer}
\begin{split}
e^{\Phi^n_{j}}_{l^2}(\kappa^j_{N+ \iota})&=\frac{|\Phi^{n,\tu{pred}}_{j}(\kappa^j_{N+\iota})-\Phi^n_{j}(\kappa^j_{N+\iota})|}{|\Phi^n_{j}(\kappa^j_{N+\iota})|}\,, \\
\end{split}
\end{equation}
where $j=1,\ldots,N_v$ and $\iota = 1,\ldots,M\,.$ This formula uses the Euclidean norm (or 2-norm) of vector defined in Section \ref{sec:model}.  Note that because each of our online multiscale basis functions is evaluated at all $m = 1089$ fine-grid nodes of the conforming fine grid $\mathcal{T}_h\,,$
the relative $l^2$ error computed by the formula \eqref{Cer} gives similar error to the relative
error in $L^2(\Omega)$-norm, that is,
\begin{align}\label{msber}
\begin{split}
e_{L^2}(\kappa^j_{N+ \iota}) = \left(\frac{\int_{\Omega} |\Phi^{n,\tu{pred}}_{j}(\kappa^j_{N+\iota})- \Phi^n_{j}(\kappa^j_{N+\iota})|^2 }{\int_{\Omega}|\Phi^n_{j}(\kappa^j_{N+\iota})|^2}\right)^{\frac{1}{2}}\,.
\end{split}
\end{align}
In this paper, for the sake of simplicity, we use the error formula \eqref{Cer} only.

As defined in Section \ref{sec:model}, for vectorization $\bfa{v}(\bfa{x})$ of a matrix, given the orthonormal basis $\psi_1, \ldots, \psi_{N_m}$, we recall the formula for vector-valued function $\bfa{v}(\bfa{x}) = (v_1, \ldots, v_{N_m}) = v_1 \psi_1 + \cdots + v_{N_m} \psi_{N_m}\,:$ 
\begin{equation}\label{normvvf}
\|\bfa{v}(\bfa{x})\|_2 = \left( \int_{\Omega}|\bfa{v}(\bfa{x})|^2 \right)^{1/2} = \left( \int_{\Omega} \sum_{l=1}^{N_m} |v_l(\bfa{x})|^2 \right)^{1/2} =\left(\sum_{l=1}^{N_m} \int_{\Omega}|v_l(\bfa{x})|^2\right)^{1/2}\,.
\end{equation}

The predicted multiscale solution is a more significant indicator of how beneficial the trained neural network is.  Thus, we compare the predicted solution $p^{n+1,\tu{pred}}_{\tu{ms}}$ obtained from \eqref{apred}--\eqref{finesolpred} with $p^{n+1}_{\tu{ms}}$ defined
by \eqref{matrixms}--\eqref{finerepsol}.  These multiscale solutions are reached at the last Picard iterative step $(n+1)$th (without the time discretization for the steady-state case \eqref{ssc}, and at the final temporal step $S$th for the time-dependent case \eqref{tdc}). The relative errors are calculated in $L^2(\Omega)$-norm and $H^1(\Omega)$-norm (both in fine-scale coordinates as \eqref{finerepsol} and \eqref{finesolpred}) based on the formula \eqref{normvvf}, that is,
\begin{align}\label{msser}
\begin{split}
e_{L^2}(\kappa^j_{N+ \iota}) &= \left(\frac{\int_{\Omega} |p^{n+1,\tu{pred}}_{\tu{ms}}(\kappa^j_{N+\iota}) - p^{n+1}_{\tu{ms}}(\kappa^j_{N+\iota})|^2 }{\int_{\Omega}|p^{n+1}_{\tu{ms}}(\kappa^j_{N+\iota})|^2}\right)^{\frac{1}{2}}\,,\\
e_{H^1}(\kappa^j_{N+ \iota}) &= \left(\frac{\int_{\Omega} |\nabla p^{n+1,\tu{pred}}_{\tu{ms}}(\kappa^j_{N+\iota}) - \nabla p^{n+1}_{\tu{ms}}(\kappa^j_{N+\iota})|^2 }{\int_{\Omega}|\nabla p^{n+1}_{\tu{ms}}(\kappa^j_{N+\iota})|^2}\right)^{\frac{1}{2}}\,.
\end{split}
\end{align}

\bigskip

Using the norms defined in \eqref{bochner} for Bochner space and \eqref{normvvf}, the relative $L^2$ and $H^1$ errors for the time-dependent problem are
\begin{align}
\label{errors_single_time}
\begin{split}
e^p_{L^2,L^2}(\kappa^j_{N+\iota})&=\frac{||p^{n+1,\tu{pred}}_{\tu{ms}}(\kappa^j_{N+\iota}) - p^{n+1}_{\tu{ms}}(\kappa^j_{N+\iota})||_{L^2(0,T;L^2(\Omega))}}{|| p^{n+1}_{\tu{ms}}(\kappa^j_{N+\iota})||_{L^2(0,T;L^2(\Omega))}}\,, \\
e^p_{L^2,H^1}(\kappa^j_{N+\iota})&=\frac{||\nabla p^{n+1,\tu{pred}}_{\tu{ms}}(\kappa^j_{N+\iota}) - \nabla p^{n+1}_{\tu{ms}}(\kappa^j_{N+\iota})||_{L^2(0,T;\bfa{L}^2(\Omega))}}{||\nabla p^{n+1}_{\tu{ms}}(\kappa^j_{N+\iota})||_{L^2(0,T;\bfa{L}^2(\Omega))}}\,.
\end{split}
\end{align}
In these two formulas, the situation is as steady-state case, only Picard iteration is employed to discretize space, and the solutions have continuous time in $[0,T]\,.$ 

We create $(6000 \cdot 81)$ realizations of heterogeneous stochastic permeability fields in experiments (thus, $6000$ distinct realizations for each of $N_v=81$ coarse neighborhoods).  Recall that Fig.~\ref{fig:kappa} presents several examples of permeability fields.  Next, $(6000 \cdot 81)$ sample pairs are generated using these $(6000 \cdot 81)$ permeability realizations.
From such $(6000 \cdot 81)$ sample pairs, we typically use $N_s \cdot N_v =5000 \cdot 81$ sample pairs as training samples (with 20\% of them designated for validation) and the remaining $M \cdot N_v =1000 \cdot 81$ sample pairs in the role of testing samples.

It should be noted that neural networks are not usually fed raw data directly.
Generally, to facilitate network optimization and improve the likelihood of successful outcomes, data must be well-prepared.  Specifically, data processing techniques as normalization and standardization are employed to rescale input and output variables prior to training a neural network model.
%
The value distribution of a standardized dataset is scaled to have a standard deviation of 1 and a mean of 0. On the other hand, normalization involves squeezing data into a provided range (typically $[0,1]$ or $[-1,1]$), which is established by the activation function being employed. 

This paper solely takes into account data normalization.  Our focus is on a somewhat straightforward normalizing procedure, namely reciprocal normalization, 
%
that normalizes values to the range of -1 to 1, by the normalization and denormalization forms
\begin{equation}
   x'= 2 \frac{x  - \tu{min}}{\tu{max} - \tu{min}} -1\,, \qquad       x = \bigg[\frac{x'+1}{2}(\text{max}-\text{min})\bigg]+\text{min}\,.
\end{equation}
The normalized value in this case is $x'$; the original, or denormalized value is $x\,;$ the maximum and minimum observable values of the provided data are $\text{max}$ and $\text{min}\,,$ respectively. 
 It should be noted that all input and target variables in this research are subject to this kind of normalization. The MinMaxScaler as scikit-learn object is employed in the normalization process of our dataset in simulations. 

In the following, we will consider two model problems and investigate how neural networks are utilized.  Also, we will test the networks' performance using online multiscale basis functions with different numbers $Nb$ of offline multiscale basis functions.


\subsection{Experiment for steady-state Richards equation}\label{ssc}

In this part, we first examine the steady-state Richards equation, which is obtained from \eqref{eq:original0}: find $p \in V$ such that
\beq
\label{eq:original00}
- \div \left(\kappa(\bfa{x}) \, \frac{1}{1+|p|}  \nabla p \right) = f \  \textrm{in } \Omega\,,
	\eeq
using the source term $f=1$ and the Dirichlet boundary condition $p(\bfa{x})=0$ over $\partial \Omega\,.$ 

With this form \eqref{eq:original00}, a dataset of $(6000 \cdot 81)$ sample pairs are generated (so 6000 different samples for each of $N_v=81$ coarse neighborhoods), by employing the terminology introduced in Section \ref{dlgms} and the online GMsFEM in Section \ref{onstage}.  Next, deep neural networks (DNNs) are constructed using $(5000 \cdot 81)$ training samples from that dataset in order to estimate the local online multiscale basis functions in the sense of prediction (instead of the traditional online GMsFEM).  For varying numbers $Nb$ of offline multiscale basis functions ($Nb= 2, 4, 6, 8, 12, 16$), Fig.~\ref{PrBatch1} plots the relative $l^2$ error RMSE \eqref{rl2e} between the DNN-predictions $\Phi^{n,\tu{pred}}_{j}$ and targets $\Phi^{n}_{j}$ in the training set (including the validation set), corresponding to different batch sizes and epoch numbers. This illustrates the convergence of the loss function \eqref{loss} for the neural networks.
\vspace{-0.3cm}
\begin{figure}[H]
\centering
\subfloat[][2 offline multiscale basis functions.]{\begin{subfigure}{0.31\textwidth}
\centering
\footnotesize{Epochs -- 100 / Batch size -- 16}\\
\includegraphics[width=0.99\linewidth]{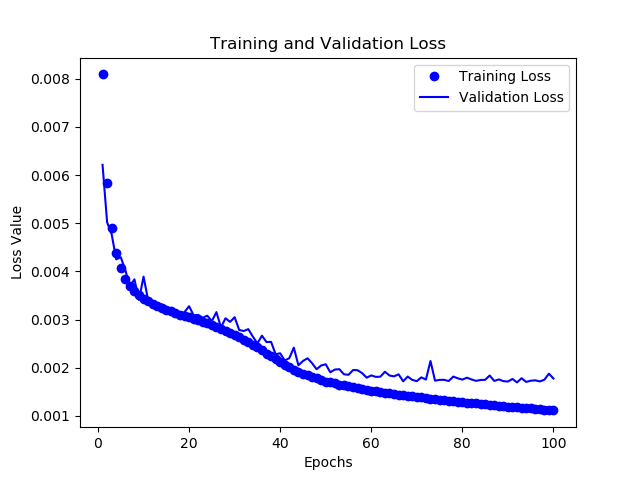}
\end{subfigure}
\begin{subfigure}{0.31\textwidth}
\centering
\footnotesize{Epochs -- 100 / Batch size -- 32}\\
\includegraphics[width=0.99\linewidth]{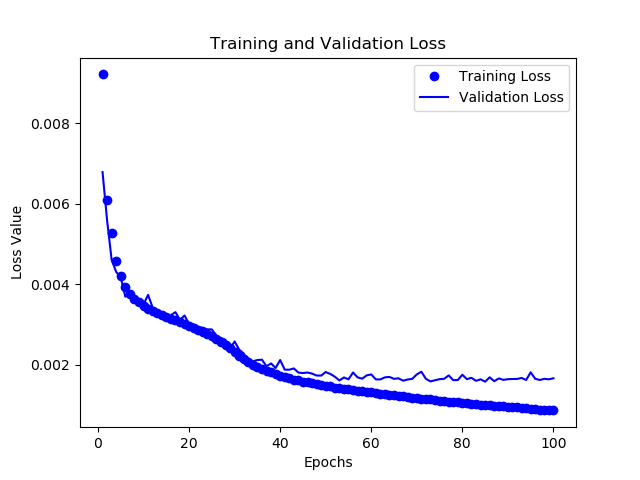}
\end{subfigure}
\begin{subfigure}{0.31\textwidth}
\centering
\footnotesize{Epochs -- 100 / Batch size -- 64}\\
\includegraphics[width=0.99\linewidth]{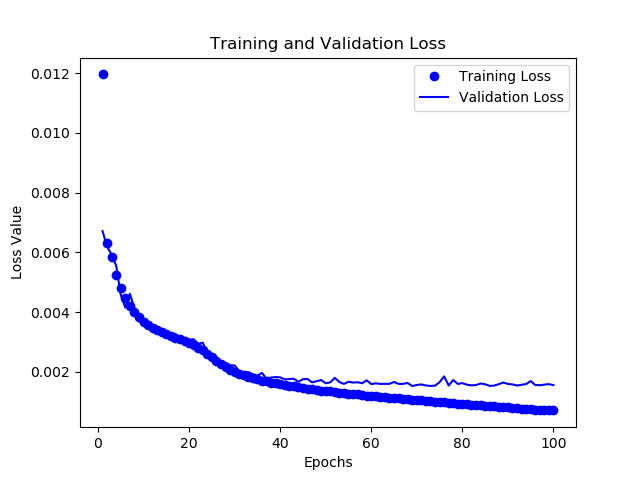}
\end{subfigure}
}
\end{figure}
\vspace{-1cm}
\begin{figure}[H]
\ContinuedFloat
\centering
\subfloat[][4 offline multiscale basis functions.]{
\begin{subfigure}{0.31\textwidth}
\centering
\footnotesize{Epochs -- 100 / Batch size -- 16}\\
\includegraphics[width=0.99\linewidth]{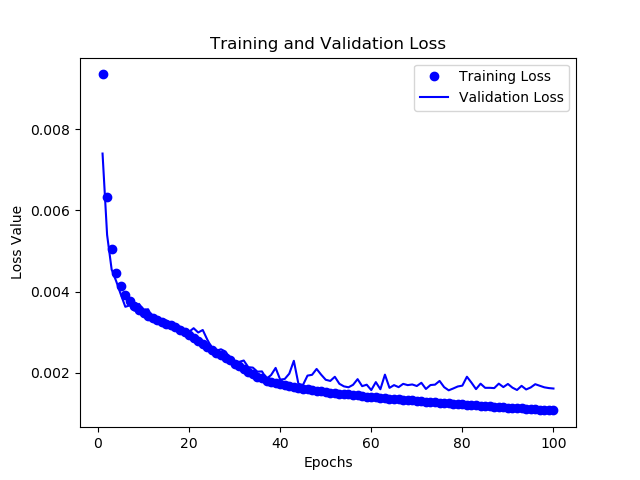}
\end{subfigure}
\begin{subfigure}{0.31\textwidth}
\centering
\footnotesize{Epochs -- 100 / Batch size -- 32}\\
\includegraphics[width=0.99\linewidth]{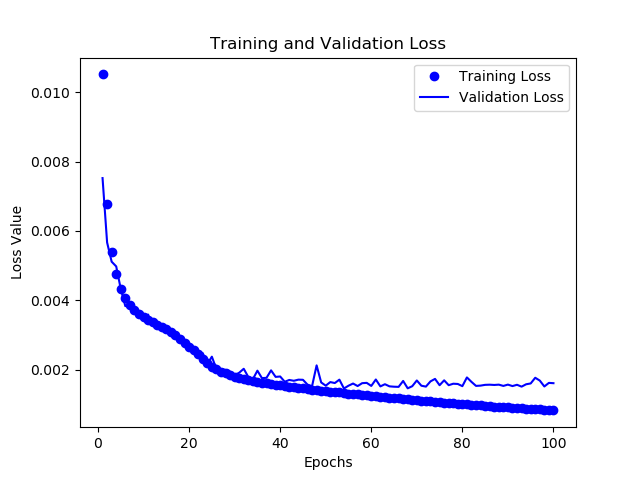}
\end{subfigure}
\begin{subfigure}{0.31\textwidth}
\centering
\footnotesize{Epochs -- 100 / Batch size -- 64}\\
\includegraphics[width=0.99\linewidth]{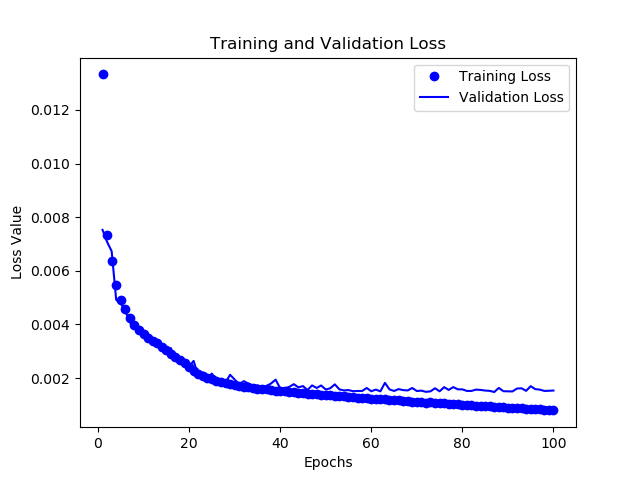}
\end{subfigure}
}
\end{figure}
\vspace{-1cm}
\begin{figure}[H]
\ContinuedFloat
\centering
\subfloat[][6 offline multiscale basis functions.]{
\begin{subfigure}{0.31\textwidth}
\centering
\footnotesize{Epochs -- 100 / Batch size -- 16}\\
\includegraphics[width=0.99\linewidth]{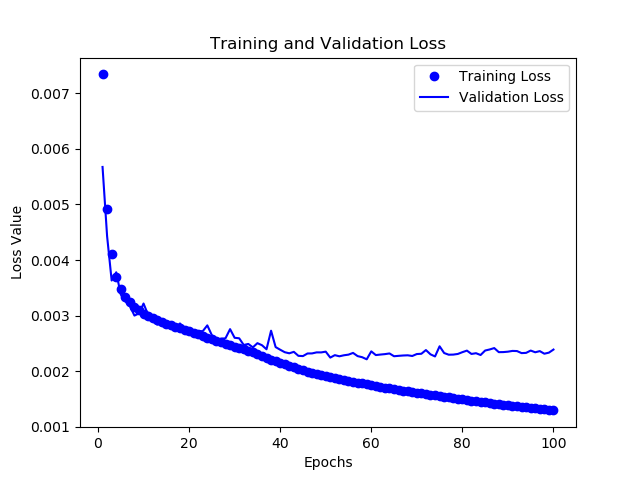}
\end{subfigure}
\begin{subfigure}{0.31\textwidth}
\centering
\footnotesize{Epochs -- 100 / Batch size -- 32}\\
\includegraphics[width=0.99\linewidth]{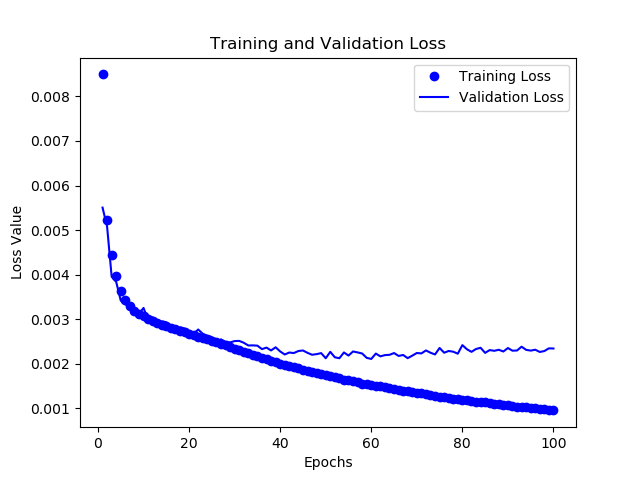}
\end{subfigure}
\begin{subfigure}{0.31\textwidth}
\centering
\footnotesize{Epochs -- 100 / Batch size -- 64}\\
\includegraphics[width=0.99\linewidth]{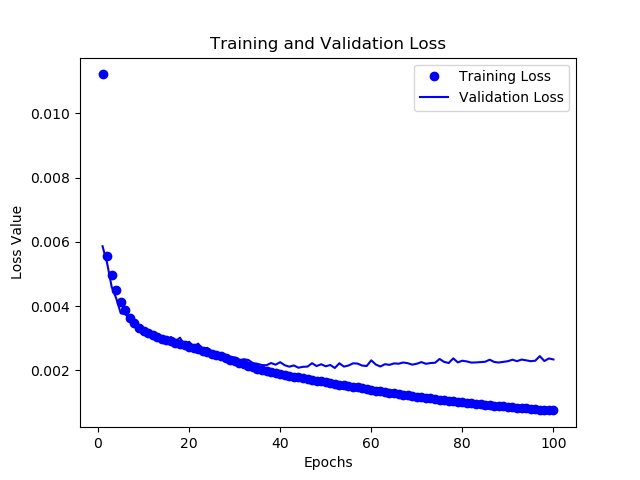}
\end{subfigure}
}
\end{figure}
\vspace{-1cm}
\begin{figure}[H]
\ContinuedFloat
\centering
\subfloat[][8 offline multiscale basis functions.]{
\begin{subfigure}{0.31\textwidth}
\centering
\footnotesize{Epochs -- 100 / Batch size -- 16}\\
\includegraphics[width=0.99\linewidth]{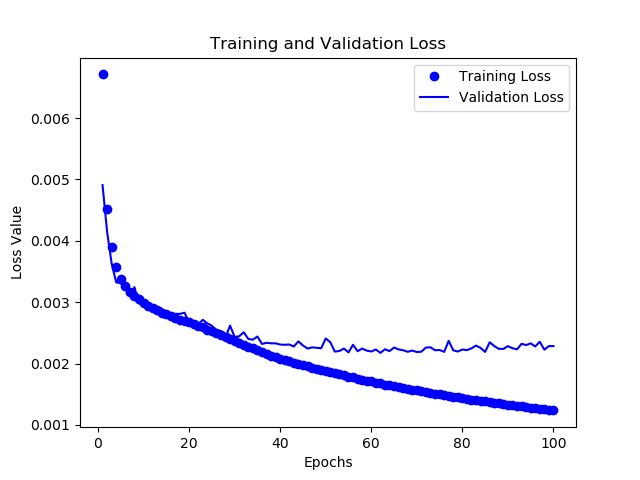}
\end{subfigure}
\begin{subfigure}{0.31\textwidth}
\centering
\footnotesize{Epochs -- 100 / Batch size -- 32}\\
\includegraphics[width=0.99\linewidth]{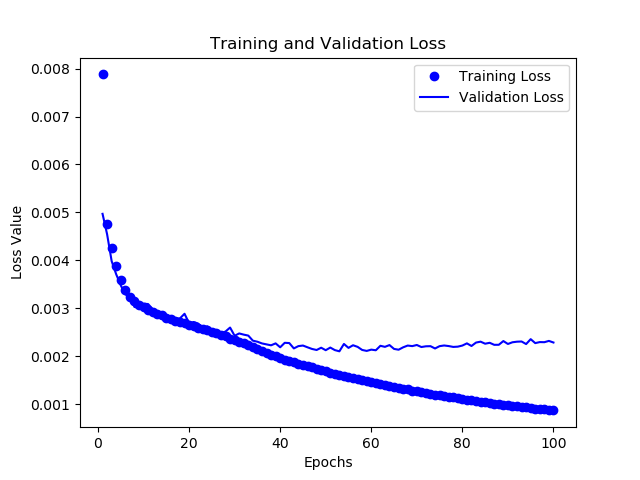}
\end{subfigure}
\begin{subfigure}{0.31\textwidth}
\centering
\footnotesize{Epochs -- 100 / Batch size -- 64}\\
\includegraphics[width=0.99\linewidth]{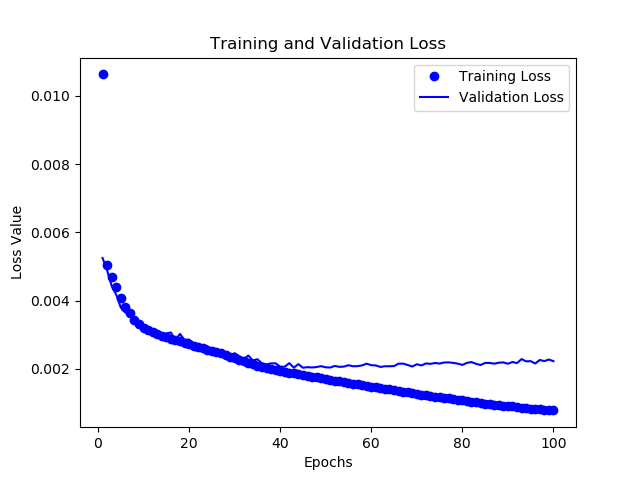}
\end{subfigure}
}
\end{figure}
\vspace{-1cm}
\begin{figure}[H]
\ContinuedFloat
\centering
\subfloat[][12 offline multiscale basis functions.]{
\begin{subfigure}{0.31\textwidth}
\centering
\footnotesize{Epochs -- 100 / Batch size -- 16}\\
\includegraphics[width=0.99\linewidth]{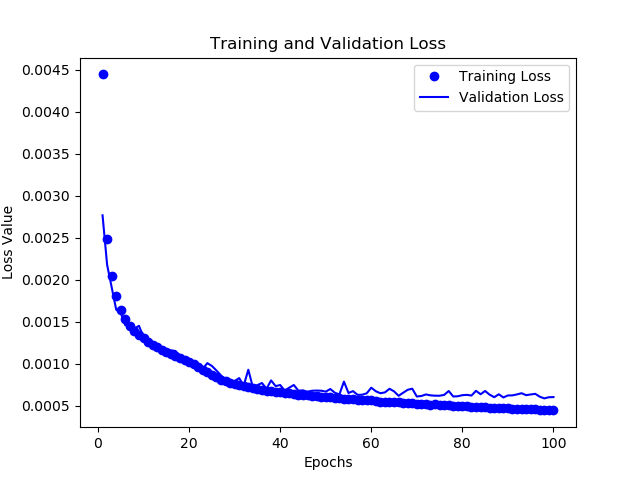}
\end{subfigure}
\begin{subfigure}{0.31\textwidth}
\centering
\footnotesize{Epochs -- 100 / Batch size -- 32}\\
\includegraphics[width=0.99\linewidth]{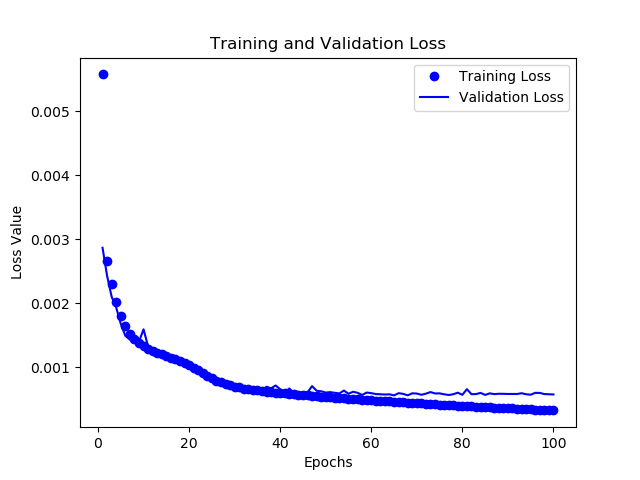}
\end{subfigure}
\begin{subfigure}{0.31\textwidth}
\centering
\footnotesize{Epochs -- 100 / Batch size -- 64}\\
\includegraphics[width=0.99\linewidth]{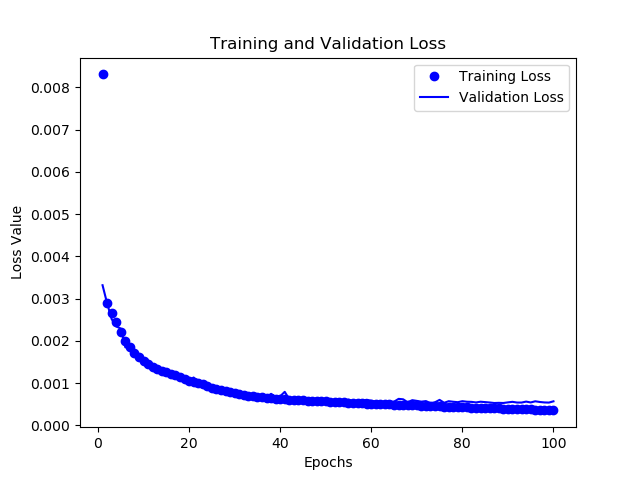}
\end{subfigure}
}
\end{figure}
\vspace{-0.8cm}
\begin{figure}[H]
\ContinuedFloat
\subfloat[][16 offline multiscale basis functions.]{
\begin{subfigure}{0.31\textwidth}
\centering
\centering
\footnotesize{Epochs -- 100 / Batch size -- 16}\\
\includegraphics[width=0.99\linewidth]{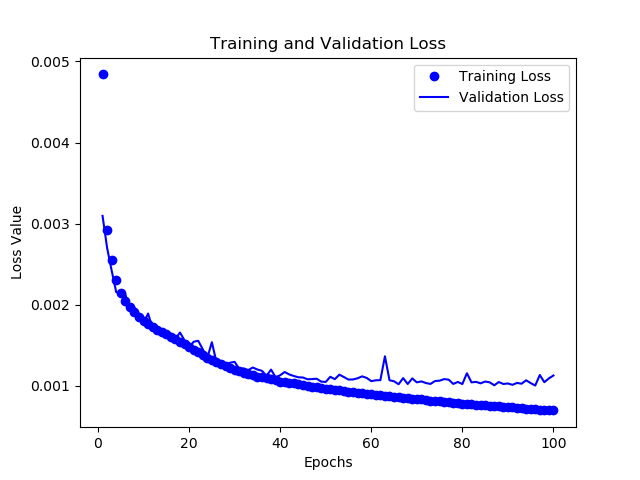}
\end{subfigure}
\begin{subfigure}{0.31\textwidth}
\centering
\footnotesize{Epochs -- 100 / Batch size -- 32}\\
\includegraphics[width=0.99\linewidth]{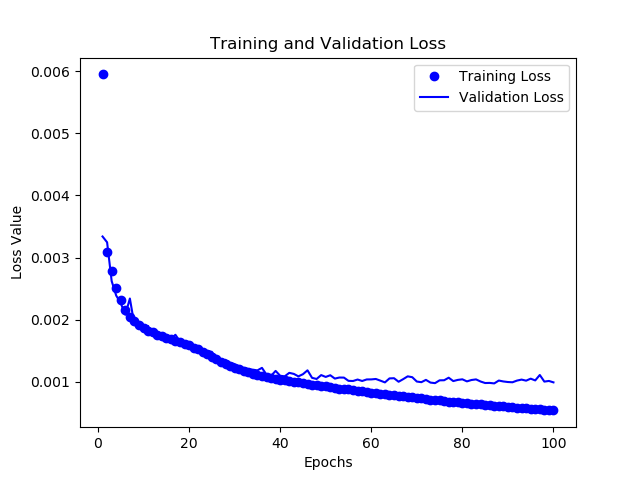}
\end{subfigure}
\begin{subfigure}{0.31\textwidth}
\centering
\footnotesize{Epochs -- 100 / Batch size -- 64}\\
\includegraphics[width=0.99\linewidth]{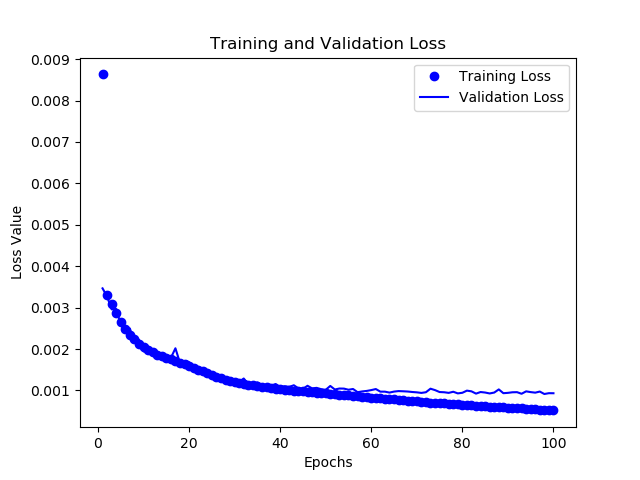}
\end{subfigure}
}
\caption{Steady-state case \eqref{eq:original00}.  Prediction's relative $l^2$ errors RMSE \eqref{rl2e} for online basis over $(5000 \cdot 81)$ training samples.}
\label{PrBatch1}
\end{figure}
\vspace{-0.8cm}
\begin{table}[H]
	\begin{center}
		\begin{tabular}{|c|c|c|c|c|c|c|c|c|c|}
			\hline
   Batch size & \multicolumn{3}{|c|}{16} & \multicolumn{3}{|c|}{32} & \multicolumn{3}{|c|}{64} \\ 
   \hline
			\backslashbox{$Nb$}{Error (\%)} &  min & max & mean & min & max & mean & min & max & mean\\
            \hline
			2 & 2.15 & 114.40 & 16.71 & 1.90 & 108.21  & 16.19 & 1.57 & 156.79  & 15.06\\
			4 & 2.19 & 64.32 & 12.03 & 1.91 & 58.50 & 11.38 & 2.04 & 52.84  & 11.61\\
			6 & 2.97 & 137.15 & 22.51 & 2.45 & 127.68 & 21.99 & 1.55 & 118.36  & 21.66\\
            8 & 2.90 & 135.11 & 22.44 & 2.37 & 152.61 & 21.19 & 2.01 & 111.73  & 20.48\\
            12 & 1.18 & 52.01 & 10.56 & 2.35 & 124.22 & 20.51 & 1.69 & 59.22  & 10.27\\
            16 & 2.04 & 108.43 & 19.80 & 1.94 & 103.96 & 17.73 & 2.07 & 179.75  & 26.35\\
			\hline
		\end{tabular}
	\end{center}
	\caption{Steady-state case \eqref{eq:original00}.  Prediction's relative $l^2$ errors \eqref{Cer} for online basis, using 100 epochs with different batch sizes, on $(1000 \cdot 81)$ testing samples.}
	\label{testres1}
\end{table}
\vspace{-0.8cm}
\begin{table}[H]
	\begin{center}
		\begin{tabular}{|c|c|c|c|c|c|c|c|c|c|}
			\hline
   Epochs & \multicolumn{3}{|c|}{50} & \multicolumn{3}{|c|}{100} & \multicolumn{3}{|c|}{150} \\ 
   \hline
			\backslashbox{$Nb$}{Error (\%)} &  min & max & mean & min & max & mean & min & max & mean\\
            \hline
			2 & 2.11 & 87.80 & 17.81 & 1.90 & 108.21  & 16.19 & 1.72 & 87.42  & 15.38\\
			4 & 1.50 & 62.74 & 11.89 & 1.91 & 58.50 & 11.38 & 1.79 & 57.99  & 11.75\\
			6 & 3.89 & 150.28 & 22.99 & 2.45 & 127.68 & 21.99 & 2.17 & 186.02  & 22.81\\
            8 & 3.30 & 191.39 & 21.31 & 2.37 & 152.61 & 21.19 & 2.36 & 188.63  & 21.39\\
            12 & 2.85 & 125.59 & 21.33 & 2.35 & 124.22 & 20.51 & 2.06 & 107.93  & 20.67\\
            16 & 2.82 & 149.12 & 50.51 & 1.94 & 103.96 & 17.73 & 2.27 & 105.33  & 27.24\\
			\hline
		\end{tabular}
	\end{center}
	\caption{Steady-state case \eqref{eq:original00}.  Prediction's relative $l^2$ errors \eqref{Cer} for online basis, using different epochs with batch size of 32, on $(1000 \cdot 81)$ testing samples.}
	\label{testres2}
\end{table}

The testing procedure is then performed for each of the numbers $Nb$ of offline multiscale basis functions (per coarse neighborhood), using the created deep neural networks and $(1000 \cdot 81)$ testing samples. Tables \ref{testres1} and \ref{testres2} contain the computation and recording of the minimum, maximum, mean, and one sample relative $l^2$ errors \eqref{Cer} for online basis, between the DNN-predictions $\Phi^{n,\tu{pred}}_{j}$ and targets $\Phi^{n}_{j}\,.$  These Tables \ref{testres1} and \ref{testres2} show that as the number of epochs increases, so does the accuracy of the predictions.  We also consider increased batch sizes, which naturally result in growing learning rates as desired.  Our neural network architecture is highly efficient in terms of accuracy.  These studies demonstrate that the best results are obtained with 100 epochs and batch size of 32, which will be used for prediction in the rest of this paper.

\begin{figure}[H]
\begin{center}
\includegraphics[width=0.62\linewidth]{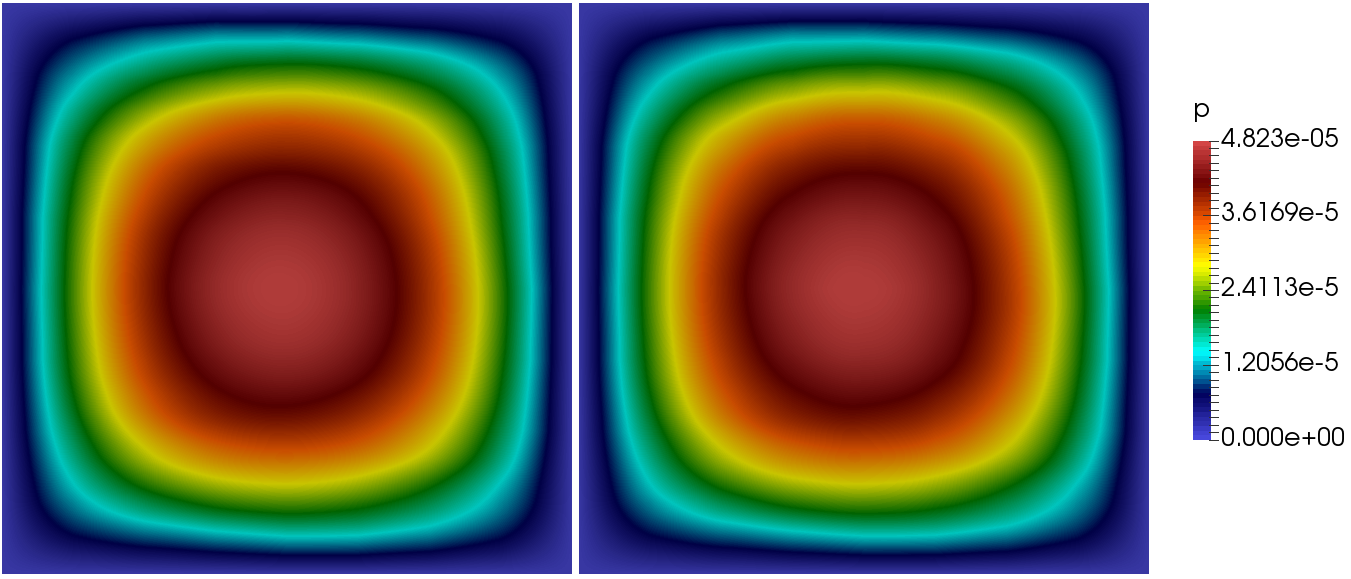}\\\hspace{-50pt} (A) $p^{n+1}_{\tu{ms}}\,.$ \hspace{150pt} (B) $p^{n+1,\tu{pred}}_{\tu{ms}}\,.$ 
	\end{center}
	\caption{Steady-state problem \eqref{eq:original00}.  Numerical results for solutions (using the first permeability field from Fig.~\ref{fig:kappa}): (A) online GMsFEM solution;
	(B) multiscale solution obtained by prediction of online multiscale basis functions. 	}
	\label{fig:steady}
\end{figure}

The multiscale solution of the steady-state problem \eqref{eq:original00} is now obtained numerically, utilizing the first permeability field in Fig.~\ref{fig:kappa}. In particular, Fig.~\ref{fig:steady} shows the solution $p^{n+1}_{\tu{ms}}$ of \eqref{finerepsol} achieved by the online GMsFEM with four offline multiscale basis functions, and the corresponding solution $p^{n+1,\tu{pred}}_{\tu{ms}}$ to \eqref{finesolpred} attained via DNN-prediction of local online multiscale basis functions.

Afterward, we complete numerical experiments by employing all $(1000\cdot 81)$ testing samples (not just one permeability field above). That is, we compare the predicted solution $p^{n+1,\tu{pred}}_{\tu{ms}}$ of \eqref{finesolpred} (obtained through the predicted online basis) with the online GMsFEM solution $p^{n+1}_{\tu{ms}}$ of \eqref{finerepsol} in relative $L^2$ and $H^1$ errors \eqref{msser}. Specifically, Table \ref{tab2} presents the comparison between $p^{n+1,\tu{pred}}_{\tu{ms}}$ of \eqref{finesolpred} and $p^{n+1}_{\tu{ms}}$ of \eqref{finerepsol} in the $L^2$ error norm, whereas Table \ref{tab3} displays the comparison in the $H^1$ error norm.  These tables show one sample errors (for the first permeability realization from Fig.~\ref{fig:kappa}) together with the mean, minimum, and maximum errors (across $1000\cdot 81$ testing samples), using different numbers $Nb$ of offline multiscale basis functions, to observe convergence of the proposed algorithm. 
The highest error is smaller than $6 \% \,,$ and all predictions are made with good precision. There is less of difference between mean errors and minimum errors than there is between mean errors and maximum errors.  In terms of accuracy, this fact shows that the majority of the numerical results are nearer to the minimum errors.  To sum up, highly precise solution to the steady-state Richards equation can be obtained by DNN-prediction.

\begin{table}[H]
	\begin{center}
		\begin{tabular}{|c|c|c|c|c|}
			\hline
			 $Nb$ & mean error (\%) & minimum error (\%) & maximum error (\%) &  one sample (\%) \\
			\hline
			2 & 0.328  & 0.189 & 0.737 & 0.307 \\
			4 & 0.144 & 0.095 & 0.275 & 0.176 \\
			6 & 0.122 & 0.084 & 0.213 & 0.115 \\
		      8 & 0.127 & 0.088 & 0.201 & 0.109 \\
                12 & 0.057 & 0.035 & 0.109 & 0.051\\
                16 & 0.098 & 0.054 & 0.185 & 0.108\\
			\hline
		\end{tabular}
	\end{center}
	\caption{Steady-state case \eqref{eq:original00}.  Numerical results for solutions over $(1000 \cdot 81)$ testing samples (the case one sample uses the first permeability field from Fig.~\ref{fig:kappa}): $L^2$ error \eqref{msser} between the online GMsFEM solution and GMsFEM solution with predicted online basis.}
	\label{tab2}
\end{table}
\vspace{-0.5cm}
\begin{table}[H]
	\begin{center}
		\begin{tabular}{|c|c|c|c|c|}
			\hline
			$Nb$ & mean error (\%) & minimum error (\%) & maximum error (\%) &  one sample (\%) \\
			\hline
			2 & 3.129 & 2.032 & 5.910 & 2.911 \\
			4 & 1.974 & 1.362 & 3.788 & 2.352 \\
			6 & 1.772 & 1.265 & 2.851 & 1.707 \\
		      8 & 1.921 & 1.355 & 3.108 & 1.689 \\
                12 & 1.304 & 0.853 & 2.112 & 1.145\\
                16 & 2.125 & 1.272 & 3.522 & 2.191\\
			\hline
		\end{tabular}
	\end{center}
	\caption{Steady-state case \eqref{eq:original00}.  Numerical results for solutions over $(1000\cdot 81)$ testing samples (the case one sample uses the first permeability field from Fig.~\ref{fig:kappa}): $H^1$ error \eqref{msser} between online GMsFEM solution and GMsFEM solution with predicted online basis.}
	\label{tab3}
\end{table}

Now, we estimate the execution time at each stage of the method (employing the first permeability field $\kappa$ from Fig.~\ref{fig:kappa}) in order to properly assess the efficiency of the proposed deep neural networks \eqref{nnbm}. 
An amount of $2004.1321$ s is required to train the neural network $\mathcal{N}_{B}^{j,n}(\kappa)$ \eqref{nnbm}. Keep in mind that the training procedure only needs to be done once, and it should not be included to our whole algorithm's execution time. Therefore, using the testing samples, it is more important to compare the time that the network $\mathcal{N}_{B}^{j,n}(\kappa)$ takes to predict one local online multiscale basis function $\Phi^{n,\tu{pred}}_{j}(\kappa)$ \eqref{bpred} with the time that the standard online GMsFEM costs to compute the corresponding local online multiscale basis function $\Phi^{n}_{j}(\kappa)$ \eqref{Phion}, that are respectively 0.023 s and 0.425744 s. Hence, the prediction is much faster than the standard for each online multiscale basis function (using one testing sample, over a coarse neighborhood). Also, selecting one testing sample, we compare the excution time for obtaining the predicted solution $p^{n+1,\tu{pred}}_{\tu{ms}}$ of \eqref{finesolpred} (containing as well the time for predicting all $N_v$ online multiscale basis functions) with the time for computing the online GMsFEM solution $p^{n+1}_{\tu{ms}}$ of \eqref{finerepsol} (including the time for generating all $N_v$ online multiscale basis functions in Section \ref{onstage}, without counting the time for creating the provided offline multiscale basis functions).  Specifically, the solution $p^{n+1,\tu{pred}}_{\tu{ms}}$ is computed in 1.295 s, whereas $p^{n+1}_{\tu{ms}}$ is reached after 2.11852 s. These results clearly indicate that the proposed method using deep neural networks provides a speedup of approximately 2 times. 
This efficiency gain is attributed to the power of our neural network $\mathcal{N}_{B}^{j,n}(\kappa)$ of rapidly generating online multiscale basis functions, and giving us the chance to solve Eq.~\eqref{eq:original00} instantly at the final Picard iterative step, in contrast to the more time-consuming process of using many Picard iterations together with the traditional online GMsFEM. Such efficiency becomes especially apparent in the cases involving a fine mesh, demonstrating the proposed DNN method's practical usefulness.  Also, it is observable that our deep neural networks can solve the problem at a much faster pace without sacrificing the correctness of the solution.  For the time-dependent problem in the next Section \ref{tdc}, this approximation of our scheme's execution speed would also hold true. 

\bigskip

\subsection{Experiment for time-dependent Richards equation}\label{tdc}

The time-dependent Richards equation \eqref{eq:original0} is investigated here: find $p \in V$ such that
\vspace{10pt}
\begin{equation} \label{time_dep_eq}
\frac{\partial p(t,\bfa{x})}{\partial t} - \div\left ( \kappa(\bfa{x}) \, \frac{1}{1+|p|}\, \nabla p(t,\bfa{x}) \right )
= f(t,\bfa{x}) \  \textrm{in} \ (0,1] \times \Omega\,,   
\end{equation}
with $f(t,\bfa{x}) = \sin(\pi x_1)\cos(\pi x_2)\,,$ having the initial condition $p(0,\bfa{x})= p_{0}=0$ in $\Omega\,,$ as well as the Dirichlet boundary condition $p(t,\bfa{x})=0$ on $(0,T] \times \partial \Omega\,.$  Recall that we chose the terminal time $T=5\cdot 10^{-5} = S\tau\,,$ and $S=20$ temporal steps, with the time step size $\tau=T/S = 25 
\cdot 10^{-7}\,.$

Using \eqref{matrixms}--\eqref{finerepsol}, we compute the local online multiscale basis functions (by the previous data from the  Picard step $n$th) and generate $(6000\cdot 81)$ samples.  Such online multiscale basis functions are approximated by deep neural networks (DNNs), which are constructed by $(5000\cdot 81)$ training samples (from the given $6000\cdot 81$ samples). Note that these $(5000\cdot 81)$  vectorized training permeability fields are input together into one neural network and need several batches.   


At the last temporal step $S$th, Fig.~\ref{fig:dynamic} presents numerical results for solving the multiscale equation derived from the time-dependent Richards equation \eqref{eq:original0}, employing four offline multiscale basis functions and  the first permeability field from Fig.~\ref{fig:kappa}. 
 That is, in Fig.~\ref{fig:dynamic}, at the final Picard iteration $(n+1)$th, the solution $p^{n+1}_{\tu{ms}}$ is obtained from \eqref{matrixms}--\eqref{finerepsol} via the online GMsFEM, whereas our DNN-prediction of online multiscale basis functions yields the solution $p^{n+1,\tu{pred}}_{\tu{ms}}$ attained from \eqref{apred}--\eqref{finesolpred}.

\begin{figure}[H]
	\begin{center}
		\begin{minipage}[h]{0.44\linewidth}
			\center{\includegraphics[width=0.8\linewidth]{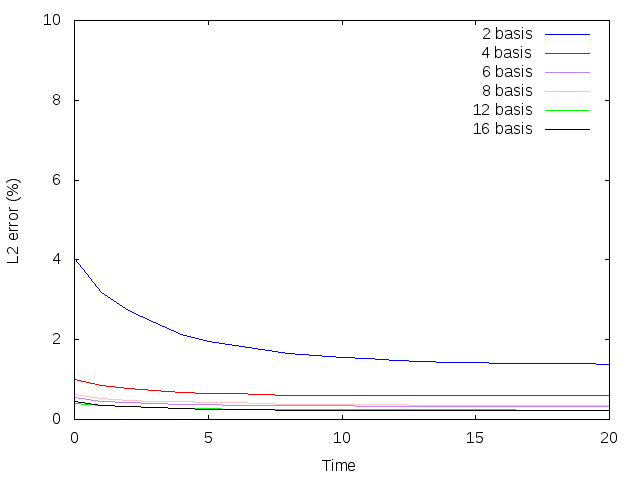}\\(A) Relative mean $L^2$ errors between $p^{n+1}_{\tu{ms}}$ and $p^{n+1,\tu{pred}}_{\tu{ms}}\,.$}
		\end{minipage}
		\begin{minipage}[h]{0.44\linewidth}
			\center{\includegraphics[width=0.8\linewidth]{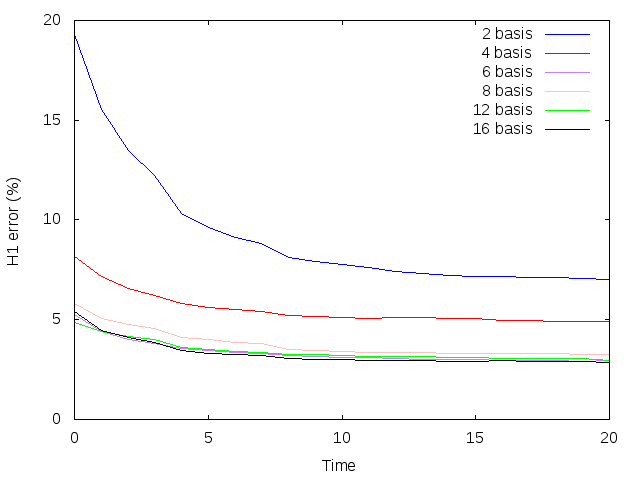}\\(B) Relative mean $H^1$ errors between $p^{n+1}_{\tu{ms}}$ and $p^{n+1,\tu{pred}}_{\tu{ms}}\,.$}
		\end{minipage}
	\end{center}
	\caption{The time-dependent problem \eqref{time_dep_eq}.  Numerical results for solutions achieved by online GMsFEM and by DNN-prediction of online basis functions, over $(1000\cdot 81)$ testing samples: (A) distribution of mean relative $L^2$ errors \eqref{msser} by time; (B) distribution of mean relative $H^1$ errors \eqref{msser} by time.}
	\label{graph1}
\end{figure}

\begin{figure}[H]
	\begin{center}
		\begin{minipage}[h]{0.44\linewidth}
			\center{\includegraphics[width=0.8\linewidth]{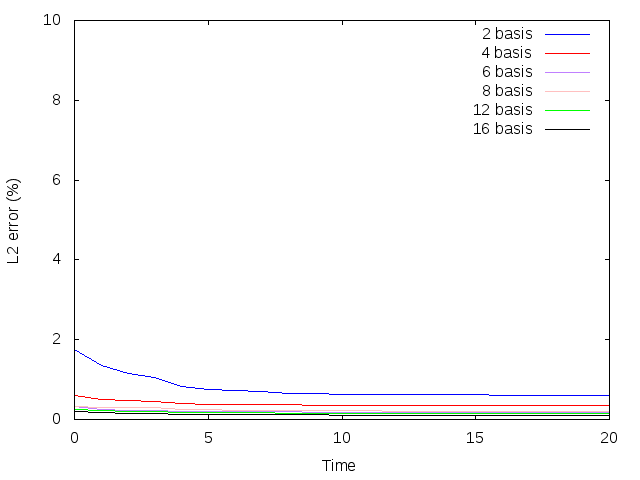}\\(A) Relative minimum $L^2$ errors between $p^{n+1}_{\tu{ms}}$ and $p^{n+1,\tu{pred}}_{\tu{ms}}\,.$}
		\end{minipage}
		\begin{minipage}[h]{0.44\linewidth}
			\center{\includegraphics[width=0.8\linewidth]{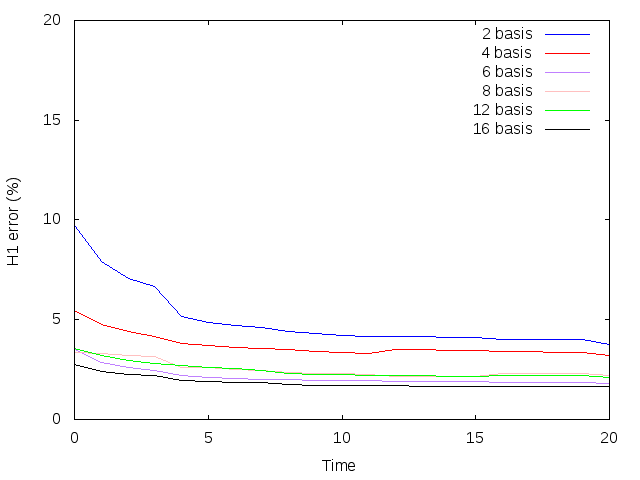}\\(B) Relative minimum $H^1$ errors between $p^{n+1}_{\tu{ms}}$ and $p^{n+1,\tu{pred}}_{\tu{ms}}\,.$}
		\end{minipage}
	\end{center}
	\caption{The time-dependent problem \eqref{time_dep_eq}.  Numerical results for solutions achieved by online GMsFEM and by DNN-prediction of online basis functions, over $(1000\cdot 81)$ testing samples: (A) distribution of minimum relative $L^2$ errors \eqref{msser} by time; (B) distribution of minimum relative $H^1$ errors \eqref{msser} by time.}
	\label{graph2}
\end{figure}

\begin{figure}[H]
	\begin{center}
		\begin{minipage}[h]{0.44\linewidth}
			\center{\includegraphics[width=0.8\linewidth]{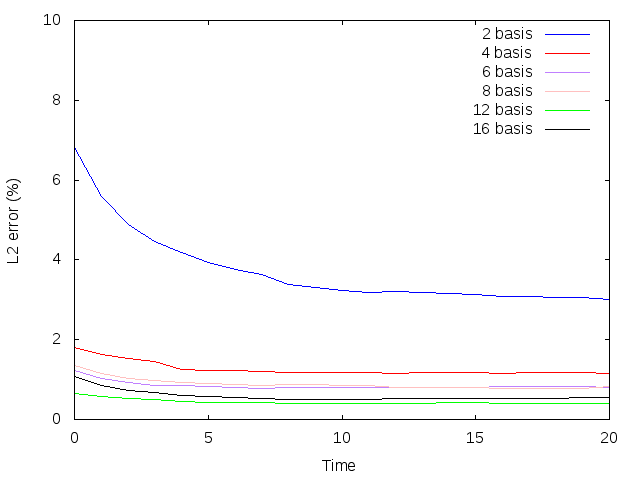}\\(A) Relative maximum $L^2$ errors between $p^{n+1}_{\tu{ms}}$ and $p^{n+1,\tu{pred}}_{\tu{ms}}\,.$}
		\end{minipage}
		\begin{minipage}[h]{0.44\linewidth}
			\center{\includegraphics[width=0.8\linewidth]{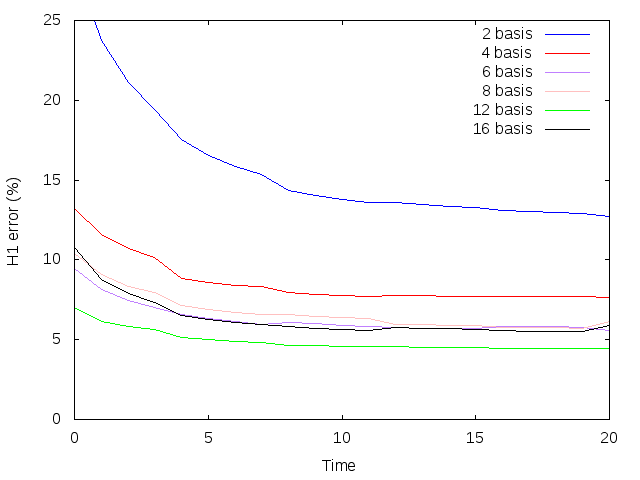}\\(B) Relative maximum $H^1$ errors between $p^{n+1}_{\tu{ms}}$ and $p^{n+1,\tu{pred}}_{\tu{ms}}\,.$}
		\end{minipage}
	\end{center}
	\caption{The time-dependent problem \eqref{time_dep_eq}.  Numerical results for solutions achieved by online GMsFEM and by DNN-prediction of online basis functions, over $(1000\cdot 81)$ testing samples: (A) distribution of maximum relative $L^2$ errors \eqref{msser} by time; (B) distribution of maximum relative $H^1$ errors \eqref{msser} by time.}
	\label{graph3}
\end{figure}

\begin{figure}[H]
	\begin{center}
		\begin{minipage}[h]{0.44\linewidth}
			\center{\includegraphics[width=0.8\linewidth]{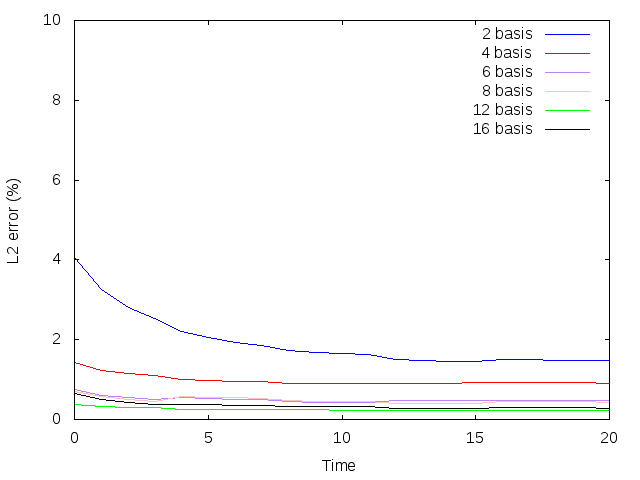}\\(A) Relative $L^2$ errors between $p^{n+1}_{\tu{ms}}$ and $p^{n+1,\tu{pred}}_{\tu{ms}}\,$ for one sample. } 
		\end{minipage}
		\begin{minipage}[h]{0.44\linewidth}
			\center{\includegraphics[width=0.8\linewidth]{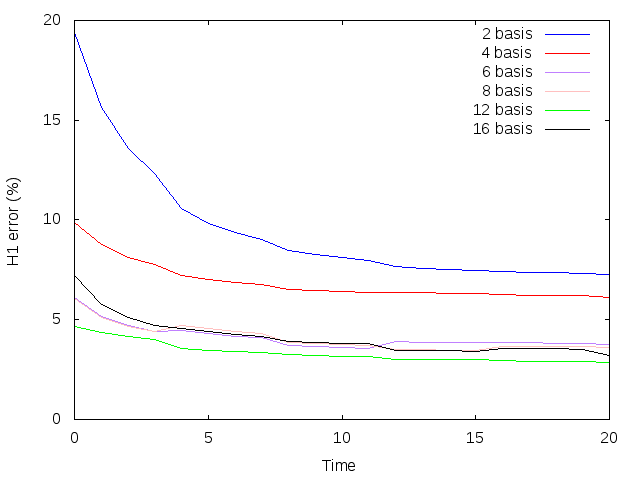}\\(B) Relative $H^1$ errors between $p^{n+1}_{\tu{ms}}$ and $p^{n+1,\tu{pred}}_{\tu{ms}}\,$ for one sample.} 
		\end{minipage}
	\end{center}
	\caption{The time-dependent problem \eqref{time_dep_eq}.  Numerical results for solutions achieved by online GMsFEM and by DNN-prediction of online basis functions, for one sample (using the first permeability field from Fig.~\ref{fig:kappa}): (A) distribution of relative $L^2$ errors \eqref{msser} by time; (B) distribution of relative $H^1$ errors \eqref{msser} by time.}
	\label{graph4}
\end{figure}

Numerical experiments are then conducted using all $(1000\cdot 81)$ testing samples for each of the developed deep neural networks corresponding to each of the numbers $Nb$ of offline multiscale basis functions.  That is, in terms of the relative $L^2$ and $H^1$ errors \eqref{msser}, we compare the online GMsFEM solution $p^{n+1}_{\tu{ms}}$ \eqref{finerepsol} with the predicted solution $p^{n+1,\tu{pred}}_{\tu{ms}}$ \eqref{finesolpred}.
For various numbers $Nb$ of offline multiscale basis functions, the relative $L^2$ and $H^1$ errors between $p^{n+1}_{\tu{ms}}$ \eqref{finerepsol} and $p^{n+1,\tu{pred}}_{\tu{ms}}$ \eqref{finesolpred} are shown in Tables \ref{tab5} and \ref{tab6} at the last time step.  Then, Tables \ref{tab7} and \ref{tab8} record such relative errors using the $L^2$ and $H^1$ norms in Bochner space \eqref{errors_single_time} for continuous time in $[0,T]\,.$ Those tables display the error distribution in mean, minimum, maximum, and one sample errors at the final Picard iteration. Also, at the testing stage, Figs.~\ref{graph1}--\ref{graph4} depict the distribution of the solution's relative $L^2$ and $H^1$ errors at the last Picard iterative step and at the selected five temporal steps (1st, 5th, 10th, 15th, 20th). It should be noted that the online multiscale basis functions (calculated directly by the online GMsFEM or predicted by DNNs) are added only at those chosen five temporal steps. Given each of these five time steps, one separate neural network is constructed for each of the numbers $Nb$ of offline basis functions.
For training, the same $(5000\cdot 81)$ samples of permeability fields are generated and are the inputs of each of these five neural networks; however, due to the five time steps, the outputs of such five networks are distinct $5\cdot (5000\cdot 81)$ local online multiscale basis functions, corresponding to different $5\cdot (5000\cdot 81)$ training samples of local online multiscale basis functions.  Similarly, for testing, the same $(1000\cdot 81)$ samples of permeability fields are the inputs of each of these five neural networks; whereas, the outputs of such five networks are different $5\cdot (1000 \cdot 81)$ local online multiscale basis functions, relating to distinct $5\cdot (1000\cdot 81)$ testing samples of local online multiscale basis functions.    

\ds{
\begin{figure}[H]
\begin{center}  \includegraphics[width=0.56\linewidth]{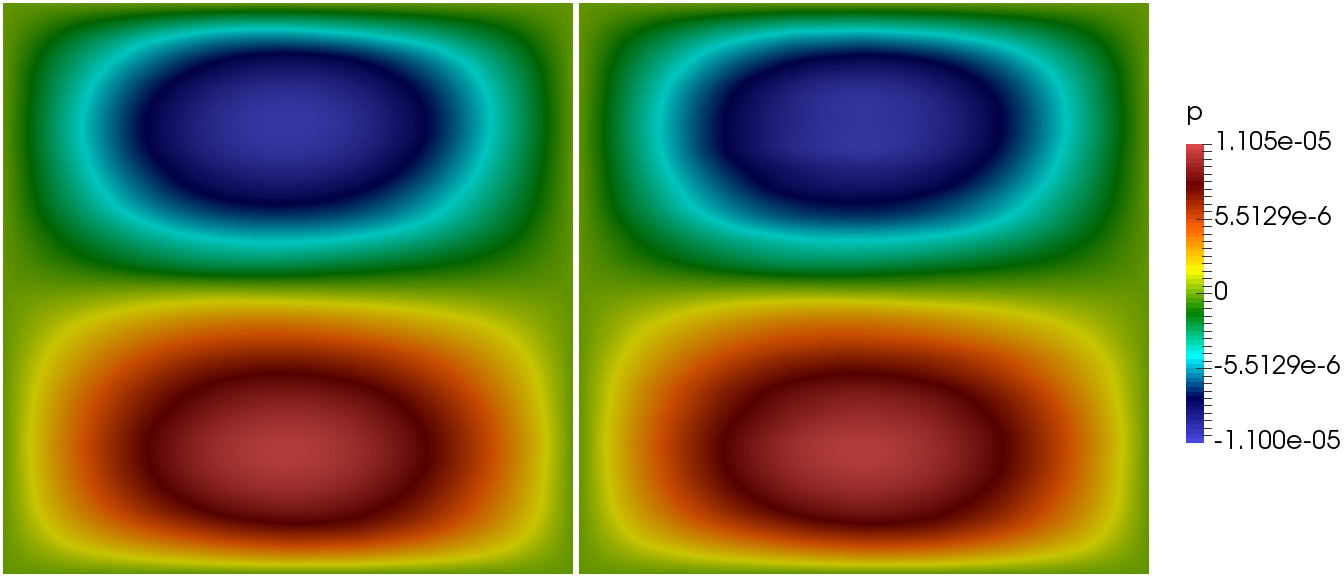}\\ \hspace{-50pt} (A) $p^{n+1}_{\tu{ms}}\,.$ \hspace{140pt} (B) $p^{n+1,\tu{pred}}_{\tu{ms}}\,.$ 
	\end{center}
	\caption{The time-dependent problem \eqref{eq:original0}.  Numerical results for solutions at the last temporal step (using the first permeability field from Fig.\ \ref{fig:kappa}): 	(A) online GMsFEM solution; (B) predicted solution.
 }
	\label{fig:dynamic}
\end{figure}}

\ds{
\begin{table}[H]
	\begin{center}
		\begin{tabular}{|c|c|c|c|c|}
			\hline
			 $Nb$ & mean error (\%) & minimum error (\%) & maximum error (\%) &  one sample (\%) \\
			\hline
			2 &  1.387 & 0.592 & 3.015 & 1.476 \\
			4 &  0.602 & 0.351 & 1.152 & 0.910 \\
			6 &  0.322 & 0.171 & 0.802 & 0.471 \\
		      8 &  0.352 & 0.203 & 0.833 & 0.432 \\
                12 & 0.231 & 0.145 & 0.407 & 0.213 \\
                16 & 0.217 & 0.106 & 0.551 & 0.265 \\
			\hline
		\end{tabular}
	\end{center}
	\caption{The time-dependent problem \eqref{time_dep_eq}.  Numerical results for solutions at the last time step over $(1000\cdot 81)$ testing samples (the case one sample uses the first permeability field from Fig.\ \ref{fig:kappa}): relative $L^2$ error \eqref{msser} between online GMsFEM solution $p^{n+1}_{\tu{ms}}$ and solution $p^{n+1,\tu{pred}}_{\tu{ms}}$ obtained by prediction of online basis functions at the last Picard step.}
	\label{tab5}
\end{table}}
\vspace{-0.5cm}
\begin{table}[H]
	\begin{center}
		\begin{tabular}{|c|c|c|c|c|}
			\hline
			$Nb$ & mean error (\%) & minimum error (\%) & maximum error (\%) &  one sample (\%) \\
			\hline
			2 & 6.995  & 3.748 & 12.728 & 7.284 \\
			4 & 4.901  & 3.214 & 7.636 & 6.121 \\
			6 & 2.941  & 1.801 & 5.563 & 3.779 \\
		      8 & 3.243  & 2.216 & 6.113 & 3.591 \\
                12 & 2.969 & 2.103 & 4.451 & 2.863 \\
                16 & 2.879 & 1.662 & 5.896 & 3.231 \\
			\hline
		\end{tabular}
	\end{center}
	\caption{The time-dependent problem \eqref{time_dep_eq}.  Numerical results for solutions at the last time step over $(1000\cdot 81)$ testing samples (the case one sample uses the first permeability field from Fig.\ \ref{fig:kappa}): relative $H^1$ error \eqref{msser} between online GMsFEM solution $p^{n+1}_{\tu{ms}}$ and solution $p^{n+1,\tu{pred}}_{\tu{ms}}$ obtained by prediction of online basis functions at the last Picard step.}
	\label{tab6}
\end{table}
\vspace{-0.5cm}
\begin{table}[H]
	\begin{center}
		\begin{tabular}{|c|c|c|c|c|}
			\hline
			 $Nb$ & mean error (\%) & minimum error (\%) & maximum error (\%) &  one sample (\%) \\
			\hline
			2 & 1.601 & 0.678 & 3.348 & 1.676 \\
			4 & 0.619  & 0.372 & 1.182 & 0.934 \\
			6 &  0.339 & 0.184 & 0.799 & 0.479 \\
		      8 & 0.377  & 0.225 & 0.836 & 0.455 \\
                12 & 0.248 & 0.158 & 0.413 & 0.233 \\
                16 & 0.231 & 0.114 & 0.536 & 0.312 \\
			\hline
		\end{tabular}
	\end{center}
	\caption{The time-dependent problem \eqref{time_dep_eq}.  Numerical results for solutions over $(1000\cdot 81)$ testing samples (the case one sample uses the first permeability field from Fig.\ \ref{fig:kappa}): relative $L^2$ error in Bochner space \eqref{errors_single_time} between online GMsFEM solution $p^{n+1}_{\tu{ms}}$ and solution $p^{n+1,\tu{pred}}_{\tu{ms}}$ obtained by prediction of online basis functions at the last Picard step. }
	\label{tab7}
\end{table}
\vspace{-0.5cm}
\begin{table}[H]
	\begin{center}
		\begin{tabular}{|c|c|c|c|c|}
			\hline
			$Nb$ & mean error (\%) & minimum error (\%) & maximum error (\%) &  one sample (\%) \\
			\hline
			2 &  8.062 & 4.495 & 14.221 & 8.321 \\
			4 & 5.189  & 3.569 & 7.841 & 6.502 \\
			6 & 3.135  & 2.006 & 5.775 & 3.901 \\
		      8 &  3.505 & 2.404 & 6.257 & 3.812 \\
                12 & 3.226 & 2.318 & 4.595 & 3.141 \\
                16 & 3.063 & 1.764 & 5.848 & 3.776 \\
			\hline
		\end{tabular}
	\end{center}
	\caption{The time-dependent problem \eqref{time_dep_eq}.  Numerical results for solutions over $(1000\cdot 81)$ testing samples (the case one sample uses the first permeability field from Fig.\ \ref{fig:kappa}): relative $H^1$ error in Bochner space \eqref{errors_single_time} between online GMsFEM solution $p^{n+1}_{\tu{ms}}$ and solution $p^{n+1,\tu{pred}}_{\tu{ms}}$ obtained by prediction of online basis functions at the last Picard step. }
	\label{tab8}
\end{table}
\vspace{-0.9cm}
\vspace{-0.9cm}

Tables \ref{tab5}--\ref{tab8} make it clear that the errors behave similarly to the steady-state problem \eqref{eq:original00}.  The only difference is that, in the time-dependent case, we observe a slight increase in error compared to the steady-state one.  Nevertheless, the errors remain at a good level and the desired accuracy can be achieved by adding more offline multiscale basis functions. For example, by employing four offline multiscale basis functions, the proposed approach already has a high level of precision. All the errors are below $8 \%$ when four or more offline multiscale basis functions are used.  However, once the total number of multiscale basis functions (both offline basis functions and online basis functions) surpasses some certain numbers, the error will not decrease anymore.  With regard to other aspect, the error values from the Tables \ref{tab7} and \ref{tab8} indicate that the resulting solution does not have any jumps in error over time and has similar error level shown in the Tables \ref{tab5} and \ref{tab6}.
This fact is observed more clearly in Figs.~\ref{graph1}--\ref{graph4}, which express the distribution of relative $L^2$ and $H^1$ errors over time. The mean error distribution, in terms of accuracy, is nearer to the minimum error than to the maximum one.  Figs.~\ref{graph1}--\ref{graph4} clearly show that adding two offline multiscale basis functions are not sufficient to obtain a highly accurate multiscale solution. It is observable that the solution converges, and to obtain a precise solution, we need at least four offline multiscale  basis  functions.  To summarize, our proposed method for the time-dependent Richards problem \eqref{time_dep_eq} offers good accuracy.

\section{Conclusions}\label{sec:conclusions}
In this work, we use deep learning in conjunction with the online generalized multiscale finite element method (online GMsFEM) to build a new coarse-scale approximation approach, for the nonlinear single-continuum Richards equation as an unsaturated flow over heterogeneous non-periodic media. An originality of this strategy is that we utilize deep neural networks (DNN) to quickly and repeatedly create local online multiscale basis functions (instead of traditionally producing them by the online GMsFEM).  More specifically, the neural networks are trained using a training set consisting of random permeability realizations and the computed corresponding local online multiscale basis functions.  Our proposed deep learning technique builds nonlinear map between those permeability fields and the local online multiscale basis functions.  That is, with this nonlinear map, the nonlinearity handling of the Richards equation is incorporated into the predicted online multiscale basis functions in a novel way, reflecting any time-dependent changes in the problem's features.  In terms of online multiscale basis function predictions and hence seeking solutions, several numerical experiments conducted on two-dimensional model problems demonstrate this approach's good performance.

\bigskip

\noindent \textbf{Acknowledgements.} 

%

The work of Spiridonov Denis receives funding from the Russian Science Foundation grant No.\ 23-71-10074, https://rscf.ru/en/project/23-71-10074/.

The research work of Sergei Stepanov is supported under the state task No.\ FSRG-2024-0003.  

Tina Mai expresses her gratitude to all the valuable support from Duy Tan University, who is going to celebrate its 30th anniversary of establishment (Nov.\ 11, 1994 -- Nov.\ 11, 2024) towards ``Integral, Sustainable and Stable Development''.

We thank Aleksandr Grigorev very much for bringing us to deep learning.


\bibliographystyle{plain}
\bibliography{r1,r2}

\end{document}